\documentclass[a4paper,10pt]{article}

\newfont{\bcb}{msbm10}
\newfont{\matb}{cmbx10}
\newfont{\got}{eufm10}

\usepackage[cp1250]{inputenc}

\usepackage{amsmath,amsthm}
\usepackage{amssymb,latexsym}
\usepackage{enumerate}
\usepackage{latexsym}
      \usepackage[all, knot]{xy}
     \xyoption{arc}

\usepackage{afterpage}
\usepackage{graphicx}

\newcommand{\matR}{\mathbb{R}}

\newcommand{\matN}{\mathbb{N}}

\DeclareMathOperator{\n}{\mathfrak{N}}

\DeclareMathOperator{\odh}{\widehat{\mathcal{Q}}}

\begin{document}
\title{On the geometric and differential properties of closed sets definable in quasianalytic structures.}
\author{Iwo Biborski}



\maketitle

\begin{abstract}
In this paper we show that the equivalences between certain properties of closed subanalytic sets proved by E. Bierstone and P. Milman in \cite{[BM-1]} hold for closed sets definable in quasianalytic o-minimal structures. In particular we prove that uniform Chevalley estimate implies a stratification by the diagram of initial exponents and further, Zariski semicontinuity of the diagram of initial exponents. We also show that the stratification by the diagram implies Zariski semicontinuity of Hilbert-Samuel function. \end{abstract}

\section{Introduction}

\vspace{2ex}
This paper is concerned with a quasi-analytic structure $\mathcal{R}$, i.e. the expansion of the real field $\mathbb{R}$ by restricted Q-analytic functions. The sets definable (with parameters) in the structure $\mathcal{R}$ are precisely those subsets of $\mathbb{R}^n$ that are globally
sub-\-quasianalytic, i.e.\ sub-\-quasianalytic in a semialgebraic compactification of $\mathbb{R}^{n}$ (sub-\-quasianalytic including infinity). Trough the paper by the definable set we mean a set definable in some fixed quasianalytic structure.

We shall investigate certain natural, metric, algebro-geometric and differential properties of closed sets definable in the structure $\mathcal{R}$, including composite function property, uniform Chevalley estimate, Zariski semicontinuity of the diagram of initial exponents, semicontinuity of the Hilbert- Samuel function, semicohernce and stratification by the diagram of initial exponents.
Our research is inspired by the famous paper \cite{[BM-1]} by E. Bierstone and P. Milman, where the equivalence of these properties is established in the classical case of subanalytic sets. In this manner, each of those properties characterizes the class of subanalytic sets that are tame from the point of view of local analytic geometry. 

We shall recall precise definitions. Fix a quasi-analytic system $\mathcal{Q}=(Q_{n})_{n\in\mathbb{N}}$ of sheaves of local $\mathbb{R}$-algebras of smooth functions on $\mathbb{R}^{n}$. For each open subset $U\subset\mathbb{R}^{n}$, $Q(U)=Q_{n}(U)$ is thus a subalgebra of the algebra $\mathcal{C}^{\infty}(U)$ of real smooth functions on $U$. By a Q-analytic function (or Q-function for short), we mean any function $f\in Q(U)$. Similarly $f=(f_{1},\dots,f_{k}):U\rightarrow\mathbb{R}^{k}$ is called Q-analytic (or a Q-map) if so are its components $f_{1},\dots,f_{k}$. The following conditions on the system $\mathcal{Q}$ are imposed:\\

1. each algebra $\mathcal{Q}(U)$ contains the restrictions of polynomials;\\

2. $\mathcal{Q}$ is closed under composition, i.e. the composition of Q-mappings is a Q-mapping, whenever it is well defined;\\

3. $\mathcal{Q}$ is closed under inverse, i.e. if $\varphi:U\rightarrow V$ is a Q-mapping between open subsets $U,V\subset\mathbb{R}^{n}$, $a\in U$, $b\in V$ and if Jacobian matrix $\frac{\partial\varphi_{i}}{\partial x_{j}}(a)$ is invertible, then there are neighborhoods $U_{a}$ and $V_{b}$ of $a$ and $b$ respectively, and Q-diffeomorphism $\psi:V_{b}\rightarrow U_{a}$ such that $\varphi\circ\psi$ is the identity mapping on $V_{b}$;\\

4. $\mathcal{Q}$ is closed under differentation;\\

5. $\mathcal{Q}$ is closed under division by a coordinate, i.e. if $f\in\mathcal{Q}(U)$ and $f(x_{1},\dots,x_{i-1},a_{i},x_{i+1},\dots,x_{n})=0$ as a function in the variables $x_{j}$, $j\neq i$, then $f(x)=(x_{i}-a_{i})g(x)$ with some $g\in\mathcal{Q}(U)$;\\

6. $\mathcal{Q}$ is quasianalytic, i.e if $f\in\mathcal{Q}(U)$ and the Taylor series $\hat{f}_{a}=0$ at $a\in U$, then $f$ is 0 in a vicinity of $a$.\\

Q-analytic maps give rise, in the ordinary manner, to the category Q of Q-manifolds and Q-maps, which is a subcategory of that of smooth manifolds and smooth maps. Similarly, Q-analytic, Q-semianalytic and Q-subanalytic sets can be defined.

Denote by $\mathcal{R}={\mathcal R}_{Q}$ the expansion of the real
field $\mathbb{R}$ by restricted Q-analytic functions, i.e.\ functions
of the form
\begin{gather*} \widetilde{f}(x) = \left\{
   \begin{array}{ll}
     f(x), & \mbox{ if } x \in [-1,1]^{n} \\
     0, & \mbox{ otherwise}
   \end{array}
   \right.
\end{gather*}
where $f(x)$ is a Q-function in the vicinity of the compact cube $[-1,1]^{n}$. The structure $\mathcal {R} = \mathcal {R}_{Q}$ is model complete and o-minimal (cf. \cite{[RSW]},\cite{[KN-1]},\cite{[KN-5]},\cite{[KN-8]}). Let us recall the definitions of properties which are studied in this paper.

\vspace{2ex}
\textbf{Composite differentiable function.}
Let $M$ and $B$ be the Q-analytic manifolds, i.e. the manifolds with Q-analytic charts. Let $\varphi:M\rightarrow N$ be a proper Q-analytic mapping.

\vspace{2ex}
\textbf{Definition 1.1.} Put
\begin{gather*}
(\varphi^{*}C^{\infty}(N))\widehat{}:=\\
\left\{f\in\mathcal{C}^{\infty}(M):\,\forall_{b\in\varphi(M)}\,\exists_{ G_{b}\in\widehat{\mathcal{Q}}_{b}}: \left(\hat{f}_{a}=\hat{\varphi}_{a}^{*}(G_{b})\,,\forall_{a\in\varphi^{-1}(b)}\right)\right\}.
\end{gather*}
We say that $f\in\mathcal{C}^{\infty}(M)$ is formally composed with $\varphi$ if $f\in(\varphi^{*}\mathcal{C}^{\infty}(N))\widehat{}$.

\vspace{2ex}
Let $\varphi:M\rightarrow\mathbb{R}^{n}$ be a proper Q-analytic mapping and let $Z\subset\mathbb{R}^{n}$ be a closed set. We denote by $\mathcal{C}^{\infty}(\mathbb{R}^{n};Z)$ the Frechet algebra of smooth functions from $\mathbb{R}^{n}$ which are flat on $Z$. Then $\varphi$ induces a homomorphism 

\begin{gather*}
\varphi^{*}:\mathcal{C}^{\infty}
(\mathbb{R}^{n};Z)\rightarrow\mathcal{C}^{\infty}(M;\varphi^{-1}(Z)).
\end{gather*} 

It is clear that $(\varphi^{*}\mathcal{C}^{\infty}(\mathbb{R}^{n}))\widehat{}\,:=(\varphi^{*}\mathcal{C}^{\infty}(\mathbb{R}^{n}))\widehat{}\, \cap\mathcal{C}^{\infty}(M;\varphi^{-1}(Z))$ is closed in $\mathcal{C}^{\infty}(M;\varphi^{-1}(Z))$.

\vspace{2ex}
\textbf{Definition 1.2.} Let $Z\subset X$ be the closed definable subsets of $\mathbb{R}^{n}$. We say that $(X,Z)$ has the composite function property if, for any proper Q-analytic mapping $\varphi:M\rightarrow\mathbb{R}^{n}$ such that $\varphi(M)=X$
\begin{gather*}
\varphi^{*}\mathcal{C}^{\infty}(\mathbb{R}^{n};Z)=(\varphi^{*}\mathcal{C}^{\infty}(\mathbb{R}^{n};Z))\widehat{}.
\end{gather*}

\vspace{2ex}
\textbf{Chevalley estimate.} Let $X\subset\mathbb{R}^{n}$ be a closed definable set. By the uniformization theorem (see Theorem 2.2) there exists a proper Q-analytic mapping $\varphi:M\rightarrow\mathbb{R}^{n}$ such that $\varphi(M)=X$. We have the following

\vspace{2ex}
\textbf{Definition 1.3.} We define the formal local ideal $\mathcal{F}_{b}(X)\subset\widehat{\mathcal{Q}}_{b}$ of $X$ at $b$ in the following way:
\begin{gather*}
\mathcal{F}_{b}(X):=\bigcap_{a\in\varphi^{-1}(b)}\mbox{Ker}\,\widehat{\varphi}_{a}^{*}.
\end{gather*}

Put
\begin{gather*}
\mu_{X,b}(f):=\sup\{p\in\mathbb{R}:|f(x)|\leq const|x-b|^{p}\,,x\in X\}\\
\nu_{X,b}(f):=\max\{l\in\mathbb{N}:f\in\mathcal{F}_{b}(X)+\widehat{m}_{b}^{l}\},
\end{gather*}
where $\hat{m}_{b}$ is the maximal ideal in $\widehat{\mathcal{Q}}_{b}$. By Chevalley Lemma (see Chapter 3), for each $b\in X$ and $k\in\mathbb{N}$, there exists $l\in\mathbb{N}$ which satisfy the following condition:\\

if $f\in\widehat{\mathcal{Q}_{b}}$ and $\mu_{X,b}(f)>l$ then $\nu_{X,b}(f)>k$.

\vspace{2ex}
\textbf{Definition 1.4.} Let $l_{X}(b,k)$ be the least $l$ as above. We call $l_{X}(b,k)$ a Chevalley estimate.

\vspace{2ex}
\textbf{The diagram of initial exponents.} We identify $\widehat{\mathcal{Q}}_{b}$ with the ring of formal power series $\mathbb{R}[[y-b]]$, where $y=(y_{1},\dots,y_{n})$.

Let $\alpha=(\alpha_{1},\dots,\alpha_{n})\in\mathbb{N}^{n}$. Then the length of $\alpha$ is a sum of its coordinates: $|\alpha|=\sum_{i=1}^{n}\alpha_{i}$. We consider an ordering on $\mathbb{N}^{n}$ defined as follows: let $\alpha,\beta\in \mathbb{N}^{n}$ and $\alpha=(\alpha_{1},\dots,\alpha_{n})$, $\beta=(\beta_{1},\dots,\beta_{n})$. Then $\alpha>\beta$ if and only if $(|\alpha|,\alpha_{1},\dots,\alpha_{n})>(|\beta|,\beta_{1},\dots,\beta_{n})$ in the lexicographic order. Let $F(y)=\sum_{\beta\in\mathbb{N}^{n}}F_{\beta}(y-b)^{\beta}$. By the support of $F$ we mean the set
\begin{gather*}
\mbox{supp}\,F:=\{\beta\in\mathbb{N}^{n}:\,F_{\beta}\neq 0\}.
\end{gather*}
We denote by $\mbox{exp}\,F$ the minimum of $\mbox{supp}\,F$ in the above ordering.

Let $I$ be an ideal in $\mathbb{R}[[y-b]]$.

\vspace{2ex}
\textbf{Definition 1.5.} We call the set
\begin{gather*}
\n(I):=\{\mbox{exp}\,F:\,F\in I\setminus\{0\}\}
\end{gather*}
the diagram of initial exponents of $I$.

\vspace{2ex}
Of course $\n(I)+\mathbb{N}^{n}=\n(I)$. Thus there exists the smallest finite set $\mathcal{B}\subset\n(I)$, called the set of vertices of $\n(I)$, such that $\n(I)=\mathcal{B}+\mathbb{N}^{n}$. We write $\mathcal{D}(n)$ for the set of all diagrams in $\mathbb{N}^{n}$.

For the ideal $I$ we define the Hilbert-Samuel functions in the following way
\begin{gather*}
H_{I}(k):=\mbox{dim}_{\mathbb{R}}\frac{\mathbb{R}[[y-b]]}{I+(y-b)^{k+1}},\,\,\,k\in\mathbb{N},
\end{gather*}
here $(y-b)$ is the maximal ideal of $\mathbb{R}[[y-b]]$. By Corollary 2.2,
\begin{gather*}
H_{I}(k)=\sharp\{\alpha\in\mathbb{N}^{n}\setminus \n(I):\,|\alpha|\leq k\}.
\end{gather*}

\vspace{2ex}
\textbf{Zariski semicontinuity.} Let $Z\subset X$ be closed definable sets in $\mathbb{R}^{n}$. Let $\Gamma$ be a prtially-ordered set.

\vspace{2ex}
\textbf{Definition 1.6.} We say that a function $\kappa: X\setminus Z\rightarrow\Gamma$ is Zariski semicontinuous relatively to $Z$, if the following two conditions hold:\\

(1) for every compact $K\subset X$, the set $\kappa((X\setminus Z)\cap K)$ is finite\\

(2) for all $\gamma\in\Gamma$, $Z_{\gamma}:=Z\cup\{b\in X\setminus Z:\kappa(b)\geq\gamma\}$ is closed definable set.

\vspace{2ex}
\textbf{Formal semicohernce.} Let $Z\subset X$ be closed definable sets in $\mathbb{R}^{n}$. We provide a definition of formal semicoherence in the similar way as \cite{[BM-1]}, Definition 1.2.

\vspace{2ex}
\textbf{Definition 1.7.} We say that $X$ is formally semicoherent relatively to $Z$, if there is a definable, locally finite stratification $X=\bigcup X_{i}$ such that $Z$ is a sum of strata and, for each stratum $X_{i}$ disjoint with $Z$, there is satisfied the following condition:\\

for each $x\in\overline{X}_{i}$, there is an open neighborhood $U$ and finitely many formal power series
\begin{gather*}
f_{ij}(\cdot,Y)=\sum_{\alpha\in\mathbb{N}^{n}}f_{ij,\alpha}(\cdot)Y^{\alpha}
\end{gather*}
such that $f_{ij,\alpha}$ are Q-analytic functions on $X_{i}\cap U$, which are definable and, for each $b\in X_{i}\cap U$, $\mathcal{F}_{b}(X)$ is generated by the formal power series
\begin{gather*}
f_{ij}(b,y-b)=\sum_{\alpha\in\mathbb{N}^{n}}f_{ij,\alpha}(b)(y-b)^{\alpha}.
\end{gather*}

\vspace{2ex}
\textbf{Definition 1.8.} We say that $X$ has a stratification by the diagram of initial exponents relatively to $Z$, if there is a locally finite definable stratification of $X$ such that $Z$ is a sum of strata and the diagram of initial exponents is constant on each stratum outside $Z$.

\vspace{2ex}
Our purpose is to investigate relations between the above properties established in Definitions 1.1-1.8. We shall prove the following quasianalytic version of Theorem 1.13 from \cite{[BM-1]}:

\vspace{2ex}
\textbf{Theorem 1.1.}\textit{ Let $X\supset Z$ be closed definable subsets of $\mathbb{R}^{n}$. Then the following properties are equivalent:}\\

\textit{(1) $(X,Z)$ has a composite function property.}\\

\textit{(2) $X$ has a uniform Chevalley estimate, i.e. for every compact $K\subset X$ there is a function $l_{K}:\mathbb{N}\rightarrow\mathbb{N}$ such that $l_{X}(b,k)\leq l_{K}(k)$ for all $b\in (X\setminus Z)\cap K$.}\\

\textit{(3) there is a definable stratification of $X$ such that $Z$ is a sum of strata and the diagram of initial exponents is constant on each stratum disjoint with $Z$.}\\

\textit{(4) The diagram of initial exponents is Zariski semicontinuous, i.e. the function $b\rightarrow\n_{b}$ is Zariski semicontinuous, where $\n_{b}=\n(\mathcal{F}_{b}(X))$ for $b\in X\setminus Z$.}\\

\textit{(5) The function $b\rightarrow H_{\mathcal{F}_{b}(X)}$ is Zariski semicontinuous relatively to $Z$.}\\

\textit{(6) $X$ is formally semicoherent relatively to $Z$.\\}

\vspace{2ex}
In this paper we present the proofs of the implications $(2)\Rightarrow (3)$, $(2)\Rightarrow (4)$ and $(3)\Rightarrow(5)$, for which it is necessary to modify the proofs by E. Bierstone and P. Milman. In comparison to the proofs from \cite{[BM-1]}, property $(4)$ cannot be directly drawn form property $(2)$. To establish the semicontinuity of the diagram of initial exponents, Birstone and Milman proved simultaneosly that certain sets $Z_{\n}$ (see Definition 4.2) are subanalytic and closed, if the uniform Chevalley estimate holds (\cite{[BM-1]}, Proposition 8.6). Their proof cannot be directly applied in the quasianalytic settings, since it relies on the fact that the ring of formal power series is faithfully flat over the ring of analytic function germs. This is unknown in the quasianalytic settings. 
In Sections 2 and 3 we recall some necessary tools, like uniformization theorem, results on stratifications and trivialization of definable map or formalism of jets. We also recall Hironaka's division algorithm and construction of generic diagram of initial exponents. 
In Section 4 we rearrange the reasoning from \cite{[BM-1]} in order to prove that uniform Chevalley estimate imply the stratification by the diagram of initial exponents. We use this result in Section 5 to prove that uniform Chevalley estimate imply semicontinuity of the diagram of initial exponents. We rely our proof on reduction to the definable arcs and fact that the ring $\mathcal{Q}_{1}$ of quasianalytic function of one variable is noetherian.

The idea of reduction to the arcs is also used in Section 6 to prove that the stratification by the diagram implies semicontinuity of Hilbert-Samuel function.

The proofs of the rest of implications can be obtained in the same way as in \cite{[BM-1]} and so they are not presented here.

\section{Necessary Tools}
\vspace{2ex}
In this section we provide some necessary tools. First of all we present the concept of covering a compact quasi-analytic set by the special cubes due to K.J. Nowak (\cite{[KN-1]}). We use it to prove uniformization theorem in general definable case. Moreover, a covering of compact quasi-analytic set enables us to prove Lemma 1, which is fundamental for further investigations and cannot be proved in the same way as in analytic settings.

Next we recall the formalism of jets (see \cite{[BM-1]}, Chapter 4) and show some useful facts about Chevalley estimate. In particular we recall Chevalley lemma and its consequences. Later on, we present stratification and trivialization theorems for the quasianalytic settings, which are adaptations of the concepts due to \L ojasiewicz presented in \cite{[Lo-1]}. We end this chapter with lemma about linear equations over the noetherian local rings.

\vspace{4ex}
\textbf{Special cubes and uniformization.} Let $M$ be a Q-manifold. Let $C\subset M$.

\vspace{2ex}
\textbf{Definition 2.1.(\cite{[KN-1]})} We call $C$ a special cube of dimension $d$ in $M$ if there exists a Q-mapping $\psi$ from the vicinity of $[-1,1]^{d}$ into $M$ such that the restriction of $\psi$ to $(-1,1)^{d}$ is a diffeomorphism onto $C$.

\vspace{2ex}
We have the following

\vspace{2ex}
\textbf{Theorem 2.1.(\cite{[KN-1]})}\textit{ A relatively compact quasi-analytic subset $S\subset M$ is a finite sum of special cubes.}

\vspace{2ex}
As an immediate corollary, we obtain the following uniformization theorem

\vspace{2ex}
\textbf{Theorem 2.2.}\textit{ Let $F\subset\mathbb{R}^{n}$ be a compact definable set. Then there exist a Q-manifold M and a Q-analytic proper mapping $\varphi:M\rightarrow\mathbb{R}^{n}$ such that $\varphi(M)=X$.}

\vspace{2ex}
\textbf{Corollary 2.1.}\textit{ Let $X\subset\mathbb{R}^{n}$ be a closed definable set. Then there exists a proper Q-analytic mapping $\varphi:M\rightarrow\mathbb{R}$ such that $X=\varphi(M)$. }
\vspace{2ex}

As a consequence of Theorem 2.2 we get the following corollary( which is a general uniformization theorem)

\vspace{2ex}
By the uniformization theorem, closed definable sets are precisely those, which are the images of Q-analytic manifolds by a proper Q-analytic map.

\vspace{2ex}
Let $\varphi:M\longrightarrow\mathbb{R}^{n}$ be a proper Q-analytic mappping. Let $\varphi(a)=b$ We write $\varphi^{\ast}_{a}$ and $\widehat{\varphi}^{\ast}_{a}$ respectively for the following ring homomorphisms:

\begin{gather*}
 \varphi^{\ast}_{a}:\mathcal{Q}_{b}\ni g\, \longrightarrow\, g\circ\varphi\in\mathcal{Q}_{a}\\ 
  \widehat{\varphi}^{\ast}_{a}:\widehat{\mathcal{Q}}_{b}\ni g\, \longrightarrow\, g\circ\widehat{\varphi}_{a}\in\widehat{\mathcal{Q}}^{\ast}_{a},
\end{gather*}
where $\widehat{\varphi}_{a}$ is a Taylor power series of mappin $\varphi$.

We make use of Theorem 2.1 also in order to prove the following fundamental

\vspace{2ex}
\textbf{Lemma 2.1.}\textit{ Let $X$ be a closed, definable set in $\mathbb{R}^{n}$ and $\varphi:M\rightarrow \mathbb{R}^{n}$ be a proper Q-analytic mapping such that $\varphi(M)=X$. Let $b\in X$ and $G\in \widehat{\mathcal{Q}}_{b}$. Put
\begin{gather*}
S:=\{a\in\varphi^{-1}(b):\, \widehat{\varphi}^{\ast}_{a}(G)=0\}.
\end{gather*}
Then $S$ is an open and closed subset of $\varphi^{-1}(b)$.}

\vspace{2ex}
\begin{proof} We rearrange the proof given by M.Birstone and P.Milman for the classical analytic case. Since we do not know whether $\widehat{\mathcal{Q}}_{b}$ is faithfully flat over $\mathcal{Q}_{b}$,  we have to use different argument for Q-analytic settings.

We can assume that $M$ is an open neighborhood $U$ of $0$ in $\mathbb{R}^{n}$, $\varphi(0)=0$ and $b=0$, $\varphi(x)=(\varphi_{1}(x),\dots,\varphi_{n}(x))$. Then
\begin{gather*}
G(\varphi(x+u)-\varphi(x))=\sum_{\beta\in\mathbb{N}^{n}}\frac{D^{\beta}G(0)}{\beta!}\left(\sum_{\alpha\in\mathbb{N}^{n}\setminus \{0\}}\frac{D^{\alpha}\varphi(x)}{\alpha!}u^{\alpha}\right)^{\beta}=\sum_{\alpha\in\mathbb{N}^{n}}\frac{H_{\alpha}(x)}{\alpha!}u^{\alpha}.
\end{gather*}

Note that each $H_{\alpha}$ is a finite linear combination of derivatives of $\varphi_{i}(x)$, $i=1,\dots,n$, thus is a Q-analytic function in some common neighborhood of 0. If $a\in\varphi^{-1}(0)$, then $\hat{\varphi}^{\ast}_{a}(G)=0$ if and only if $H_{\alpha}(a)=0$ for all $\alpha$. Thus $S$ is a closed subset as an intersection of $\bigcap_{\alpha\in\mathbb{N}^{m}} H_{\alpha}^{-1}(0)$ with $\varphi^{-1}(0)$.

Observe that $G(y+v)-G(y)\in (v)\cdot\mathbb{R}[[y,v]]$, where $(v)=(v_{1},\dots,v_{n})$ is an ideal generated by $v_{i}$, $i=1,\dots,n$. Then
\begin{gather*}
G(\varphi(x+u)-\varphi(x))-G(\varphi(x+u))\in(\varphi(x))\cdot\mathbb{R}[[x,u]],
\end{gather*}
where $(\varphi(x))$ is the ideal generated by $\varphi_{i}(x)$, $i=1,\dots,n$. Suppose that $a=0\in S$, which means $\hat{\varphi}^{\ast}_{0}(G)=0$. Since $G\in \mbox{Ker}\, \widehat{\varphi}^{\ast}_{0}$, thus $G(\varphi(x+u))=0$ and $G(\varphi(x+u)-\varphi(x))\in(\varphi(x))\cdot\mathbb{R}[[x,u]]$. Thus we have

\begin{gather*}
G(\varphi(x+u)-\varphi(x))=\sum_{\alpha\in\mathbb{N}^{n}}\frac{H_{\alpha}(x)}{\alpha!}u^{\alpha}=f_{1}(x,u)\varphi_{1}(x)+\dots+f_{n}(x,u)\varphi_{n}(x),
\end{gather*}
for some $f_{i}(x,u)\in\mathbb{R}[[x,u]]$. Therefore each $H_{\alpha}$ is of the form
\begin{gather*}
H_{\alpha}(x)=\sum_{j=1}^{n}\hat{f_{j}}(x)\varphi_{i}(x),
\end{gather*}
for some $\hat{f_{i}}(x)\in\mathbb{R}[[x]]$. 

Since $\varphi$ is a proper Q-analytic map then $\varphi^{-1}(0)$ is a compact subset of $M$. By Theorem 2.1, we can represent $\varphi^{-1}(0)$ as a finite sum of special cubes:
\begin{gather*}
\varphi^{-1}(0)=\bigcup_{i}S_{i},
\end{gather*}
where, for each $i$, $S_{i}$ is a special cube, $\psi_{i}:V_{i}\rightarrow M$ is a Q-mapping from vicinity $V_{i}$ of $[-1,1]^{d_{i}}$ and $\psi|_{(-1,1)^{d_{i}}}$ is a diffeomorphism onto $S_{i}$. Since each $[-1,1]^{d_{i}}$ can be covered by intervals, therefore we can cover $\varphi^{-1}(0)$ by Q-analytic arcs, via mapping $\psi_{i}$.

Let $U$ be the neighborhood of 0 in $\mathbb{R}^{m}$. We consider quasianalytic arcs which cover $\varphi^{-1}(0)\cap U$. Let $\gamma:(-\epsilon,\epsilon)\rightarrow\mathbb{R}^{m}$ be a Q-analytic arc, $\gamma\subset\varphi^{-1}(0)\cap U$. Then
\begin{gather*}
H_{\alpha}(\gamma(t))=\sum_{i=1}^{n}\hat{f_{i}}(\gamma(t))\varphi_{i}(\gamma(t)).
\end{gather*}
Since $\gamma\subset\varphi^{-1}(0)$, $\varphi_{i}(\gamma(t))=0$, and thus $H_{\alpha}$ vanishes along $\gamma$. Therefore $H_{\alpha}$ vanishes in $\varphi^{-1}(0)\cap U$ for all $\alpha\in\mathbb{N}^{m}$, and thus $S$ is an open set in $\varphi^{-1}(0)$.
\end{proof}

Lemma 2.1 is essential for the study of the behavior of the ideal of formal power series vanishing on the germ of closed definable set $X$. In particular this ideal is a finite intersection of kernels of appropriate homomorphisms. It is included in the corollary below

\vspace{2ex}
\textbf{Corollary 2.2.}\textit{ Let $X$ be a closed definable subset of $\mathbb{R}^{n}$. Let $\varphi:M\rightarrow\mathbb{R}^{n}$ be a proper Q-analytic mapping such that $\varphi(M)=X$. Let $b\in X$ and let $s$ be the number of connected components of $\varphi^{-1}(b)$. Then}
\begin{gather*}
\mathcal{F}_{b}(X)=\bigcap_{i=1}^{s}\mbox{Ker}\,\widehat{\varphi}_{a_{i}}^{*},
\end{gather*}
\textit{where each $a_{i}$ is in a distinct connected component of $\varphi^{-1}(b)$.}

\vspace{4ex}
\textbf{Formalism of jets.}
We recall the formalism of jets in the quasianalytic settings. This construction is similar to the formalism of jets provided in \cite{[BM-1]}. 

Let $N$ be a Q-analytic manifold. Let $b\in N$, $l\in\mathbb{N}$. We write $J^{l}(b)$ for $\widehat{\mathcal{Q}}_{N,b}/\widehat{m}_{N,b}^{l+1}$, and for $G\in\widehat{\mathcal{Q}}_{N,b}$, we denote by $J^{l}G(b)$ a class of $G$ in $J^{l}(b)$.
Let $\varphi:M\rightarrow N$ be a Q-analytic mapping and let $b\in\varphi(M)\subset N$. Then for any $a\in\varphi^{-1}(b)$, the homomorphism $\widehat{\varphi}^{*}_{a}:\widehat{\mathcal{Q}}_{N,b}\rightarrow\widehat{\mathcal{Q}}_{M,a}$ induces a linear transformation 

\begin{gather*}
J^{l}\varphi(a):J^{l}(b)\rightarrow J^{l}(a).
\end{gather*}

Let $N=\mathbb{R}^{n}$ and $y=(y_{1},\dots,y_{n})$ be the coordinates in $\mathbb{R}^{n}$. We can identify $\widehat{\mathcal{Q}}_{N,b}$ with the ring of formal power series $\mathbb{R}[[y-b]]$ and thus we can treat $J^{l}(b)$ as $\mathbb{R}^{p}$, where $p={n+l\choose l}$, and $J^{l}G(b)=(D^{\beta}G(b))_{|\beta|\leq l}$.
Put 

\begin{gather*}
J^{l}_{b}:=J^{l}(b)\otimes_{\mathbb{R}}\widehat{\mathcal{Q}}_{N,b}=\bigoplus_{|\beta|\leq l}\mathbb{R}[[y-b]],
\end{gather*} 
and let $G(y)=\sum \frac{G_{\beta}}{\beta!}(y-b)^{\beta}$. We write $J^{l}_{b}G(y)$ for $(D^{\beta}G(y))_{|\beta|\leq l}\in J^{l}_{b}$. Let us notice that $J^{l}G(b)$ is a vector of constant terms in $J^{l}_{b}G(y)$.

Consider a Q-analytic mapping $\varphi:M\rightarrow \mathbb{R}^{n}$. Take $a\in M$, and let $\varphi(a)=b$. If $x=(x_{1},\dots,x_{m})$ is a system of coordinates on $M$ in a neighborhood of $a$, we identify $J^{l}(a)$ with $\mathbb{R}^{p}$, $p={m+l\choose l}$. Then
\begin{gather*}
(*)\,\,\,\,\,J^{l}\varphi(a):(G_{\beta})_{|\beta|\leq l}\longmapsto\left(\sum_{|\beta|\leq|\alpha|}G_{\beta}L_{\alpha}^{\beta}(a)\right)_{|\alpha|\leq l},
\end{gather*}
where $L^{\beta}_{\alpha}(a)=\frac{1}{\beta!}\frac{\partial^{\alpha}\varphi^{\beta}}{\partial x^{\alpha}}(a)$. By $\varphi^{\beta}$ we mean $\varphi_{1}^{\beta_{1}}\cdot\, \dots\, \cdot\varphi_{n}^{\beta_{n}}$, where $\varphi=(\varphi_{1},\dots,\varphi^{n})$ and $\beta=(\beta_{1},\dots,\beta_{n})$.

By the chain rule, a homomorphism $\widehat{\varphi}_{a}^{*}:\widehat{\mathcal{Q}}_{b}\rightarrow\widehat{\mathcal{Q}}_{a}$ induces a homomorphism $J^{l}(b)\rightarrow J^{l}(a)$ over the ring homomorphism $\widehat{\varphi}^{*}_{a}$, and thus an $\mathbb{R}[[x-a]]$-homomorphism
\begin{gather*}
J^{l}_{a}\varphi:J^{l}(b)\otimes_{\mathbb{R}}\widehat{\mathcal{Q}}_{a}\rightarrow J^{l}(a)\otimes_{\mathbb{R}}\widehat{\mathcal{Q}}_{a},
\end{gather*}
where, for any $G\in\widehat{\mathcal{Q}}_{b}$, $J^{l}_{a}\varphi((\widehat{\varphi}_{a}^{*}(D^{\beta}G))_{|\beta|\leq l})=(D^{\alpha}(\widehat{\varphi}^{*}_{a}(G)))_{|\alpha|\leq l}$. Let us notice, that we identify $J^{l}(b)\otimes_{\mathbb{R}}\odh_{a}$ with $\bigoplus_{|\beta|\leq l}\mathbb{R}[[x-a]]$ and $J^{l}(a)\otimes_{\mathbb{R}}\odh_{a}$ with $\bigoplus_{|\alpha|\leq l}\mathbb{R}[[x-a]]$. By evaluating at $a$, $J^{l}_{a}\varphi$ induces $J^{l}\varphi(a):J^{l}(b)\rightarrow J^{l}(a)$. We can identify $J^{l}\varphi(a)$ with a matrix, which coefficients are the entries of the Taylor series of $D^{\alpha}\varphi^{\beta}/\beta!$ at $a$, for $|\alpha|\leq l$ and $|\beta|\leq l$.

Let $M$ be a Q-analytic manifold, $\varphi:M\rightarrow\mathbb{R}^{n}$ a Q-analytic map, and $s\in\mathbb{N}\setminus\{0\}$. Set
\begin{gather*}
M^{s}_{\varphi}:=\{\underline{a}=(a_{1},\dots,a_{s})\in M^{s}:\varphi(a_{1})=\dots=\varphi(a_{s})\}.
\end{gather*}
We call $M^{s}_{\varphi}$ an s-fold fibre product of $M$ with respect to $\varphi$. Notice that $M^{s}_{\varphi}$ is a closed Q-analytic subset of $M^{s}$. We write $\underline{\varphi}^{s}$ for a natural map $\underline{\varphi}^{s}:M^{s}_{\varphi}\rightarrow\mathbb{R}^{n}$ such that $\underline{\varphi^{s}}(\underline{a})=\varphi(a_{1})$. Thus $\underline{\varphi}=\underline{\varphi^{s}}=\varphi\circ\pi_{i}$, where $\pi_{i}:M^{s}_{\varphi}\rightarrow M$ is the canonical projection.

In Chapter 4 we shall need the following homomorphism
\begin{gather*}
J^{l}_{\underline{a}}\varphi:J^{l}(b)\otimes_{\mathbb{R}}\odh_{M^{s}_{\varphi,\underline{a}}}\rightarrow\bigoplus_{i=1}^{s}J^{l}(a_{i}) \otimes_{\mathbb{R}}\odh_{M^{s}_{\varphi,\underline{a}}},
\end{gather*}
It is the homomorphism with $s$ components:
\begin{gather*}
J^{l}(b)\otimes\widehat{\mathcal{Q}}_{M^{s}_{\varphi},\underline{a}}\rightarrow J^{l}(a_{i})\otimes\widehat{\mathcal{Q}}_{M^{s}_{\varphi},\underline{a}},\,\,i=1,\dots,s,
\end{gather*}
each of which is obtained from the homomorphism
\begin{gather*}
J^{l}_{a_{i}}\varphi:J^{l}(b)\otimes\widehat{\mathcal{Q}}_{a_{i}}\rightarrow J^{l}(a_{i})\otimes\widehat{\mathcal{Q}}_{a_{i}}
\end{gather*}
by change of base $\widehat{\pi}^{*}_{i}:\widehat{\mathcal{Q}}_{a_{i}}\rightarrow\widehat{\mathcal{Q}}_{M^{s}_{\varphi},\underline{a}}$.
We identify $J^{l}_{\underline{a}}\varphi$ with
\begin{gather*}
J^{l}_{\underline{a}}\varphi:\bigoplus_{|\beta|\leq l}\odh_{M^{s}_{\varphi,\underline{a}}}\rightarrow\bigoplus_{i=1}^{s}\bigoplus_{|\alpha|\leq l} \otimes_{\mathbb{R}}\odh_{M^{s}_{\varphi,\underline{a}}}.
\end{gather*}

Let $L$ be the germ at $\underline{a}$ of Q-analytic subspace of $M^{s}_{\varphi}$. Then we have a homomorphism
\begin{gather*}
J^{l}_{\underline{a}}\varphi:J^{l}(b)\otimes_{\mathbb{R}}\odh_{L,\underline{a}}\rightarrow\bigoplus_{i=1}^{s}J^{l}(a_{i}) \otimes_{\mathbb{R}}\odh_{L,\underline{a}}
\end{gather*}
By evaluation at $\underline{a}$ we obtain
\begin{gather*}
J^{l}\varphi(\underline{a}):J^{l}(b)\rightarrow\bigoplus_{i=1}^{s}J^{l}_{a_{i}}.
\end{gather*}

\vspace{4ex}
\textbf{Chevalley estimate.}
We recall some notation from \cite{[BM-1]}, Chapter 5. Let $M$ be a Q-analytic manifold. Consider $M_{\varphi}^{s}$, where $\varphi:M\rightarrow\mathbb{R}^{n}$ is a proper Q-analytic mapping and $s\in\mathbb{N}$. For $\underline{a}=(a_{1},\dots,a_{s})\in M_{\varphi}^{s}$ we define(see \cite{[BM-1]}, p. 748) an ideal
\begin{gather*}
\mathcal{R}_{\underline{a}}:=\bigcap_{i=1}^{s}\mathcal{R}_{a_{i}}=\bigcap_{i=1}^{s}\mbox{Ker}\,\widehat{\varphi}_{a_{i}}^{*}.
\end{gather*}
For $k\in\mathbb{N}$, put
\begin{gather*}
\mathcal{R}^{k}(\underline{a}):=\frac{\mathcal{R}_{\underline{a}}+\widehat{\underline{m}}^{k+1}_{\underline{\varphi}(\underline{a})}}{\widehat{m}^{k+1}_{\underline{\varphi}(\underline{a})}}.
\end{gather*}
Let $b\in\mathbb{R}^{n}$ and $\Pi^{k}(b):\widehat{\mathcal{Q}}_{b}\rightarrow J^{k}(b)$ be the canonical projection. If $l\in\mathbb{N}$ and $l\geq k$, then we write $\Pi^{lk}(b):J^{l}(b)\rightarrow J^{k}(b)$ for the canonical projection of jets. For the linear transformation $J^{l}\varphi(\underline{a}):J^{l}(\underline{\varphi}(\underline{a}))\rightarrow\bigoplus_{i=1}^{s}J^{l}(a_{i})$ we write
\begin{gather*}
E^{l}:=\mbox{Ker}\,J^{l}\varphi(\underline{a}),\\
E^{lk}:=\Pi^{lk}(\underline{\varphi}(\underline{a}))E^{l}(\underline{a}).
\end{gather*}
We have the following

\vspace{2ex}
\textbf{Lemma 2.2.(Chevalley's Lemma,\cite{[BM-1]}, Lemma 5.2)}\textit{ Let $\underline{a}=(a_{1},\dots,a_{s})\in M_{\varphi}^{s}$. For all $k\in\mathbb{N}$, there exists $l\in\mathbb{N}$ such that $\mathcal{R}^{k}(\underline{a})=E^{lk}(\underline{a})$, or equivalently such that if $G\in\widehat{\mathcal{Q}}_{\underline{\varphi}(\underline{a})}$ and $\widehat{\varphi}_{a_{i}}^{*}(G)\in\widehat{m}^{l+1}_{a_{i}}$ for $i=1,\dots,s$, then $G\in\mathcal{R}_{\underline{a}}+\widehat{m}^{k+1}_{\underline{\varphi}(\underline{a})}$. }

\vspace{2ex}
Let $y=(y_{1},\dots,y_{n})\in\mathbb{R}^{n}$. Let $X$ be a closed definable set in $\mathbb{R}^{n}$.

\vspace{2ex}
\textbf{Definition 2.2.} Let $F\in\widehat{\mathcal{Q}}_{b}$. We define
\begin{gather*}
\mu_{X,b}(F):=\mbox{sup}\,\{p\in\mathbb{R}:\,|T_{b}^{l}F(y)|\leq \mbox{const}\,|y-b|^{p},\\y\in X,\,l=\min\{k\in\mathbb{N},k\geq p\}\}\\
\nu_{X,b}(F):=\max\{l\in\mathbb{N}:\,F\in\widehat{m}_{b}^{l}+\mathcal{F}_{b}(X)\}.
\end{gather*}

\vspace{2ex}
\textbf{Remark 2.1.}(\cite{[BM-1]}, Remark 6.3) It is true that $\mu_{X,b}(F)\leq\nu_{X,b}(F)$.

\vspace{2ex}
\textbf{Definition 2.3.} Let $\varphi:M\rightarrow \mathbb{R}^{n}$ be a proper Q-analytic mapping, $\varphi(M)=X$. Put
\begin{gather*}
l_{X}(b,k):=\min\{l\in\mathbb{N}: (F\in\widehat{\mathcal{Q}}_{b},\,\mu_{X,b}(F)>l)\Rightarrow\nu_{X,b}(F)>k\},\\
l_{\varphi^{*}}(b,k):=\min\{l\in\mathbb{N}:[F\in\widehat{\mathcal{Q}}_{b},\forall_{a\in\varphi^{-1}(b)}\nu_{M,a}(\widehat{\varphi}^{*}_{a}(F))>l]\Rightarrow \nu_{X,b}(F)>k\}.
\end{gather*}

\vspace{2ex}
If $\underline{a}=(a_{1},\dots,a_{s})$, $a_{i}\in\varphi^{-1}(b)$, then we define
\begin{gather*}
l_{\varphi^{*}}(\underline{a},k):=\min\{l\in\mathbb{N}:[F\in\widehat{\mathcal{Q}}_{b},\forall_{i=1\dots s}\,\,\,\nu_{M,a_{i}}(\widehat{\varphi}^{*}_{a_{i}}(F))>l]\Rightarrow\\ F\in\bigcap_{i=1}^{s}\mbox{Ker}\,\widehat{\varphi}_{a_{i}}^{*}+\widehat{m}_{b}^{k+1}\}.
\end{gather*}

\vspace{2ex}
\textbf{Remark 2.2.} By Chevalley Lemma, $l_{\varphi^{*}}(\underline{a},k)<\infty$. If $\underline{a}$ is an s-tuple of elements $a_{i}$ such that each $a_{i}$ lies in a distinct connected component of $\varphi^{-1}(b)$, thus $\bigcap_{i=1}^{s}\mbox{Ker}\,\widehat{\varphi}_{a_{i}}^{*}=\mathcal{F}_{b}(X)$, and $l_{\varphi^{*}}(b,k)\leq l_{\varphi^{*}}(\underline{a},k)$.

\vspace{2ex}

The lemma below is a quasianalytic version of Lemma 6.5 from \cite{[BM-1]}.

\vspace{2ex}
\textbf{Lemma 2.3.}\textit{ Let $\varphi:M\rightarrow\mathbb{R}^{n}$ be a proper Q-analytic mapping such that $\varphi(M)=X$. Then we have:}

(1)$\,l_{X}(b,\cdot)\leq l_{\varphi^{*}}(b,\cdot),$

(2)$for\,every\,compact\,K\subset X,\, there\,exists\,r\geq 1,\\such\,that\, l_{\varphi^{*}}(b,\cdot)\leq r\,l_{X}(b,\cdot),b\in K$.

\vspace{2ex}
In order to proof Lemma 2.3 it is enough to repeat the proof of Lemma 6.5 from \cite{[BM-1]}, since following quasianalytic version of \L ojasiewicz inequality with parameter from \cite{[KN-2]} holds:

\vspace{2ex}
\textbf{Lemma 2.4.}(\L ojasiewicz inequality with parameter,\cite{[KN-2]})\textit{ Let $f,g:A\rightarrow\mathbb{R}$ be the functions defined on a definable set $A\subset\mathbb{R}^{n}_{u}\times\mathbb{R}^{m}_{x}$. Assume that $A_{u}:=\{x\in\mathbb{R}^{n}:(u,x)\in A\}$ is a compact set for each $u\in\mathbb{R}^{n}$ and}
\begin{gather*}
f_{u},g_{u}:A_{u}\rightarrow\mathbb{R}^{n},\,\,\,f_{u}(x):=f(u,x),\,\,\,g_{u}(x):=g(u,x),
\end{gather*}
\textit{are continuous. If $\{f=0\}\subset\{g=0\}$, then there exist an exponent $\lambda>0$ and definable function $c:\mathbb{R}^{n}\rightarrow(0,\infty)$ such that}
\begin{gather*}
|f(u,x)|\geq c(u)|g(u,x)|^{\lambda},\,\,\,(u,x)\in A
\end{gather*}
\textit{for some $i=1,\dots,k$.}

\vspace{4ex}
\textbf{Stratification and trivialization.} Here we recall some trivialization and stratification theorems given by \L ojasiewicz (\cite{[Lo-1]}). Let us notice that this results hold in o-minimal structures with smooth cell decomposition. The original ideas by Hardt (\cite{[Ha]}), which relies on Weierstrass preparation theorem, cannot be adapted to the quasianalytic settings. Finally, we recall several corollaries, which for the classical analytic case were formulated in \cite{[BM-1]}. Stratification and trivialization results are based on \L ojasiewcz theorems on equitriangulations. His ideas rely on good direction and so called f-property of the set (see \cite{[Lo-1]}).

\vspace{2ex}
\textbf{Theorem 2.4.}\textit{(Trivialization)}\textit{ Let $f:E\rightarrow \mathbb{R}^{n}$ be a definable, continuous map with bounded graph. There exists partition $\mathcal{T}$ of $f(E)$ into definable leaves, such that for any $T\in\mathcal{T}$ there is a definable homeomorphism $h:f^{-1}(T)\rightarrow F\times T$, $F$--definable subset of $\mathbb{R}^{n}$, such that the following diagram is commutative:}

\begin{gather*}
\xymatrix{
f^{-1}(T) \ar[rd]_f \ar[r]^h\,\,&
T\times F \ar[d]^{\pi} \\&T}
\end{gather*}
\textit{where $\pi$ is the canonical projection.}

\vspace{2ex}

Let $M$ and $N$ be Q-analytic manifolds. Let $\Pi:M\times N\rightarrow N$.

\vspace{2ex}
\textbf{Lemma 2.5.}\textit{ Let $E\subset M\times N$ be a definable set which is relatively compact in order to $N$ and, for each $t\in M$, $E_{t}$ is finite. Then there exists a locally finite partition $\Gamma$ of $E$ into definable leaves such that, for any $T\in\Gamma$, $\Pi(T)$ is a definable leaf and $\Pi|_{T}$ is a Q-analytic isomorphism.}

\vspace{2ex}
\textbf{Theorem 2.5.}\textit{ (Stratification of a definable mappings.)}\textit{ Let $f:E\rightarrow N$ be definable, continuous and proper mapping over a closed definable set $E\subset M$. Let $\mathcal{F}$ and $\mathcal{G}$ be the locally finite families in $M$ and $N$ respectively. Then there exists a definable stratifications $\mathcal{A}$ and $\mathcal{B}$ of $M$ and $N$, compatible with $\mathcal{F}$ and $\mathcal{G}$, such that for each $A\in\mathcal{A}$ such that $A\subset E$, $f(A)\in\mathcal{B}$ and $f|_{A}:A\rightarrow f(A)$ is a trivial submersion.}

\vspace{2ex}
Theorem 2.4 and 2.5 lead to the quasianalytic version of Theorem 7.1 from \cite{[BM-1]}:

\vspace{2ex}
\textbf{Corollary 2.4.}\textit{(Stratified trivialization of definable mappings.)}\textit{ Let $E$ be a closed definable subset of Q-analytic manifold $M$. Let $\varphi:E\rightarrow\mathbb{R}^{n}$ be a proper, continuous definable mapping and let $X=\varphi(E)$. Let $\mathcal{F}$ and $\mathcal{G}$ be finite families of $E$ and $X$ respectively. Then there exist definable stratifications $\mathcal{S}$ and $\mathcal{T}$ of $E$ and $X$, such that }

\textit{(1) For each $S\mathcal\in{S}$, $\varphi(S)\in\mathcal{T}$,}

\textit{(2) For each $T\in\mathcal{T}$ and $b\in T$, there is a definable stratification $\mathcal{P}$ of $\varphi^{-1}(b)$ and a definable homeomorphism $h$ such that te following diagram is commutative:}

\begin{gather*}
\xymatrix{
\varphi^{-1}(T) \ar[rd]_\varphi \ar[r]^h\,\,&
T\times\varphi^{-1}(b)\ar[d]^{\pi} \\&T}
\end{gather*}
\textit{and for each $S\subset\varphi^{-1}(T)$, $h|_{S}$ is a Q-analytic isomorphism onto $T\times P$ for some $P\in\mathcal{P}$,}

\textit{(3) $\mathcal{S}$ is compatible with $\mathcal{F}$ and $\mathcal{T}$ is compatible with $\mathcal{G}$.}

\vspace{2ex}
\textbf{Corollary 2.5.}\textit{ If $\varphi$ is stratified as in Corollary 2.4, then, for each stratum $T$, the number of connected components of the fibre $\varphi^{-1}(b)$, where $b\in T$, is constant on $T$.}

\vspace{2ex}
Corollary 2.5 is an immediate consequence of Corollary 2.4.

Let us consider a closed definable set $E\subset M$ and a continuous definable map $\varphi:E\rightarrow\mathbb{R}^{n}$. Then, for $s\in\mathbb{N}$, the subset $E_{\varphi}^{s}$ is a closed definable subset of s-fold fibre product $M^{s}$ and we can consider a definable mapping $\underline{\varphi}:E_{\varphi}^{s}\rightarrow\mathbb{R}^{n}$. We have the following

\vspace{2ex}
\textbf{Corollary 2.6.}\textit{ Let $(\mathcal{S},\mathcal{T})$ be a stratification of $\varphi$ as in Corollary 2.4. Take $T\in\mathcal{T}$ and denote as $\mathcal{S}_{T}^{s}$ the family of all nonempty sets of the form $(S_{1}\times\dots\times S_{s})\cap E_{\varphi}^{s}$, for $S_{i}\in\mathcal{S}$ such that $S_{i}\subset\varphi^{-1}(T)$. Then $\mathcal{S}_{T}^{s}$ is a definable stratification of $\underline{\varphi}^{-1}(T)$ such that each $\underline{S}\in\mathcal{S}_{\varphi}^{s}$ admits a Q-analytic isomorphism $h:\underline{S}\rightarrow T\times P$ commuting with projection on $T$, where $P$ is a bounded definable leaf in $M^{s}$.}

Let $\underline{E}_{\varphi}^{s}$ denote a set of those $x=(x_{1},\dots,x_{s})\in E_{\varphi}^{s}$ such that each $x_{i}$ is in a distinct connected component of $\underline{\varphi}^{-1}(\underline{\varphi}(\underline{x}))$. From Corollaries 2.5 and 2.6, we immediately obtain the following

\vspace{2ex}
\textbf{Corollary 2.7.}\textit{ Let $\varphi:E\rightarrow\mathbb{R}^{n}$ be a proper definable
 map defined on a closed definable set $E\subset M$. Then $\underline{E}_{\varphi}^{s}$ is a definable subset of $M^{s}$.}

\vspace{4ex}
\textbf{ Linear equations over noetherian local rings.} Let $\mathbb{R}[[y]]$ be a ring of the formal power series, $y=(y_{1},\dots,y_{n})$ and let $(A,m)$ be a local ring such that $\widehat{A}=\mathbb{R}[[y]]$ and $\mathbb{R}[[y]]$ is faithfully flat over $A$. Let $\Phi$ be a matrix with coefficients in $A$. Let $p$ and $q$ be the number of columns and rows( respectively) of $\Phi$.  We have the following

\vspace{2ex}
\textbf{Lemma 2.6.}\textit{ Let $\widehat{\xi}\in (\mathbb{R}[[y]])^{p}$ and $\Phi\cdot\widehat{\xi}=0$. Then there exists $\xi\in A^{p}$ such that $\Psi\cdot\xi=0$ and $\widehat{\xi}-\xi\in \widehat{m}$.}
\begin{proof}
Let $\Phi_{1},\dots,\Phi_{q}$ be the rows of $\Phi$. Since $\Phi\cdot\widehat{\xi}=0$ thus $\Phi_{i}\cdot\widehat{\xi}=0$ for each $i=1,\dots,q$. There exist $\widehat{a}\in (\widehat{m})^{p}$, where $\widehat{m}$ is a maximal ideal $(y)$ in $\mathbb{R}[[y]]$, and $u\in\mathbb{R}^{p}$ such that $\Phi_{i}\cdot\widehat{a}=-\Phi_{i}\cdot u$. By the faithful flatness of $\mathbb{R}[[y]]$ over $A$  we have (see \cite{[Bo]}, Commutative Algebra, Ch. I,\S 3,5 , prop. 10 (ii)):
\begin{gather*}
\left[\begin{array}{c}\Phi_{1}\\ \Phi_{2}\\ \vdots \\ \Phi_{q} \end{array}\right]\cdot \widehat{m}^{p}\cap A^{p}=\left[\begin{array}{c}\Phi_{1}\\ \Phi_{2}\\ \vdots \\ \Phi_{q} \end{array}\right]\cdot m^{p}.
\end{gather*}
Thus there exists $a\in (m)^{p}$ such that $\Phi_{i}\cdot a=-\Phi_{i}\cdot u$. Put $\xi=a+u\in A^{p}$. Therefore $\Phi\cdot\xi=0$ as desired.
\end{proof}

\section{Generic diagram of initial exponents.}
In this chapter we recall for the reader's convenience the division algorithm of Grauert--Hironaka and the construction of the generic diagram of initial exponents from \cite{[BM-1]}, which carries over verbatim to the quasianalytic settings.

\vspace{4ex}
\textbf{ Division algorithm in the ring of formal power series.} Let $\mathbb{R}[[y]]$ be the ring of formal power series, where $y=(y_{1},\dots,y_{n})$. We have the following

\vspace{2ex}
\textbf{Theorem 3.1.}\textit{(Hironaka's division algorithm, \cite{[BM-1]}, Theorem 3.1.)}\textit{ Let $F_{1},\dots,F_{s}\in\mathbb{R}[[y]]\setminus\{0\}$ and $\alpha^{i}=\mbox{exp}\,F_{i}$, $i=1,\dots, s$. Let $\Delta_{i}:=(\alpha^{i}+\mathbb{N}^{n})\setminus \bigcup_{j=1}^{i}\Delta_{j}$ and $\Delta:=\mathbb{N}^{n}\setminus \bigcup_{j=1}^{s}\Delta_{s}$. For any $G\in\mathbb{R}[[y]]$, there exist unique $Q_{i}$, $i=1,\dots,s$, $R\in\mathbb{R}[[y]]$ such that $G=\sum_{i=1}^{s}F_{i}Q_{i}+R$, $\mbox{supp}\,Q_{i}\in\Delta_{i}$ and $\mbox{supp}\,R\in\Delta$. Moreover, $\mbox{exp}\,R\geq\mbox{exp}\,G$, $\alpha^{i}+\mbox{exp}\,Q_{i}\geq \mbox{exp}G$.}

\vspace{2ex}
\textbf{Corollary 3.1.}\textit{(\cite{[BM-1]}, Corollary 3.2)}\textit{ Let $\mathcal{I}$ be an ideal in $\mathbb{R}[[y]]$ and $\alpha^{1},\dots,\alpha^{s}$ be the vertices of $\n(\mathcal{I})$. Let $F_{1},\dots,F_{s}\in\mathcal{I}$ such that $\mbox{exp}\,F_{i}=\alpha^{i}$. Let $\Delta_{i}$ and $\Delta$ be the sets as in last theorem. Then}
\begin{gather*}
(1)\,\n(\mathcal{I})=\bigcup_{i=1}^{s}\Delta_{i},\\
(2)\,\text{There is a unique set of generators}\, G_{1},\dots,G_{s}\in \mathcal{I}\text{ such that }\\ \mbox{supp}\,(G_{i}-y^{\alpha^{i}})\subset\Delta_{i}.
\end{gather*}
The system $G_{1},\dots,G_{s}$ is called a standard basis of $I$. From Corollary 3.1, $\mathbb{R}[[y]]=I\oplus\mathbb{R}^{\n}$, where $\n$ is a diagram of initial exponents of $I$ and $\mathbb{R}[[y]]^{\n}:=\{F\in\mathbb{R}[[y]]:\,\mbox{supp}\,F\subset\mathbb{N}^{n}\setminus\n\}$.

We also recall a corollary which indicates the connection between the diagram of initial exponents and the Hilbert-Samuel function. We have the following

\vspace{2ex}
\textbf{Corollary 3.2.}\textit{(\cite{[BM-1]}, Corollary 3.2)}\textit{ Let $H_{I}$ be a Hilbert-Samuel function of $\mathbb{R}[[y]]/I$. Then, for every $k\in\mathbb{N}$}
\begin{gather*}
H_{I}(k):=\sharp\{\beta\in\mathbb{N}^{n}:\,\beta\notin\n(I)\,and\,|\beta|\leq k\}.
\end{gather*}
\textit{It follows that $H_{I}(k)$ coincides with polynomial in $k$ for sufficiently large $k$.}

\vspace{4ex}
\textbf{Lemma from linear algebra. }We recall here an useful lemma from \cite{[BM-1]}, Chapter 2. Let $V$ and $W$ be the modules over a commutative ring $R$.

\vspace{2ex}
\textbf{Definition 3.1.} Let $B\in\mbox{Hom}_{R}\,(V,W)$ and $r\in\mathbb{N}$. We define

\begin{gather*}
\mbox{ad}^{r} B\in \mbox{Hom}_{R}\left(W,\mbox{Hom}_{R}({\bigwedge}^{r}\,V,{\bigwedge}^{r+1}\,W)\right)
\end{gather*}
by formula

\begin{gather*}
(\mbox{ad}^{r}B)(\omega)(\eta_{1}\wedge\dots\wedge\eta_{r})=\omega\wedge B\eta_{1}\wedge\dots\wedge B\eta_{r},
\end{gather*}
where $\omega\in W$ and $\eta_{1},\dots,\eta_{r}\in V$.

\vspace{2ex}
\textbf{Remark 3.1.} It is clear that if $r>\mbox{rank}\, B$, then $\mbox{ad}^{r} B=0$ and if $r=\mbox{rank}\,B$ then $\mbox{ad}^{r}B\cdot B=0$.

\vspace{2ex}
We recall the following:

\vspace{2ex}
\textbf{Lemma 3.1.}\textit{(\cite{[BM-1]}, Lemma 2.1)}\textit{ Let $V$ and $W$ be finite-dimensional vector spaces over a field $\mathbb{K}$. Let $B:V\rightarrow W$ be a linear transformation. Let $r:=\mbox{rank}\, B$. Then}
\begin{gather*}
\mbox{Im}\, B=\mbox{Ker}\,\mbox{ad}^{r}B.
\end{gather*}
\textit{If $A$ is a linear transformation with target $W$, then $A\xi+B\eta=0$ if and only if} $\xi\in\mbox{Ker}\,\mbox{ad}^{r}B\cdot A$.

\vspace{2ex}
\textbf{The generic diagram of initial exponents.} In order to introduce generic diagram of initial exponents and generic Hilbert-Samuel function, we formulate the following

\vspace{2ex}
\textbf{Lemma 3.2.}\textit{ Let $A$ be a matrix of dimension $k\times n$ and $B$ be a matrix of dimension $k\times m$. Consider a block matrix $(A,B)$. Let $\pi$ be a projection from $\mathbb{R}^{n+m}$ onto $\mathbb{R}^{n}$. Then $\pi(\mbox{Ker}\,(A,B))=\{x\in\mathbb{R}^{n}:\,Ax\in \mbox{Im}\, B\}$.}

\vspace{2ex}
 \begin{proof} Let $x\in\mathbb{R}^{n}$. Then $x\in\pi(\mbox{Ker}\,(A,B))$ if and only if, there exists $y\in\mathbb{R}^{m}$ such that $(x,y)\in\mbox{Ker}\,(A,B)$. Hence $Ax+By=0$, and thus $Ax=-By$. Therefore
\begin{gather*}
x\in\pi(\mbox{Ker}\,(A,B))\Leftrightarrow \exists y\in\mathbb{R}^{m}:Ax=-By\Leftrightarrow Ax\in\mbox{Im}\,B.
\end{gather*}
\end{proof}

Let $L\subset M^{q}_{\varphi}$ be a Q-analytic leaf. Let  $\xi=(\xi_{\beta})_{|\beta|\leq l}$. We can write it as $\xi=(\xi^{k},\zeta^{lk})$, where $\xi^{k}:=(\xi_{\beta})_{|\beta|\leq k}$ and $\zeta^{lk}:=(\xi_{\beta})_{k<|\beta|\leq l}$. According to this decomposition we can write $J^{l}\varphi(\underline{a})$ as a block matrix
\begin{gather*}
J^{l}\varphi(\underline{a})=(S^{lk}(\underline{a}),T^{lk}(\underline{a}))=\left [\begin{matrix}J^{k}(\underline{a})&0\\ \star&\star\end{matrix}\right ],
\end{gather*}
for $\underline{a}\in L$ (see Chapter 1, \textbf{Formalism of jets}). By Lemma 3.2, we have
\begin{gather*}
E^{lk}(\underline{a})=\{\xi^{k}=(\xi_{\beta})_{|\beta|\leq k}:\, S^{lk}(\underline{a})\cdot\xi^{k}\in \mbox{Im}\,T^{lk}(\underline{a})\}.
\end{gather*}
Therefore $E^{lk}(\underline{a})=\mbox{Ker}\,\Theta^{lk}(\underline{a})$, $d^{lk}=\mbox{rk}\,\Theta^{lk}(\underline{a})$, where
\begin{gather*}
\Theta^{lk}(\underline{a}):=\mbox{ad}^{r^{lk}}T^{lk}(\underline{a})\cdot S^{lk}(\underline{a})\\
r^{lk}(\underline{a}):=\mbox{rk}\,T^{lk}(\underline{a}).
\end{gather*}
Summing up, we obtain the following

\vspace{2ex}
\textbf{Lemma 3.3}\textit{ Let $k\in\mathbb{N}$. For sufficiently large $l\in\mathbb{N}$ there is equality} 

\begin{gather*}
\mbox{dim}\,E^{lk}(\underline{a})=\mbox{dim}\left ( \mathcal{R}^{k}(\underline{a})+\widehat{m}_{\underline{\varphi}(\underline{a})}\right )/\widehat{m}_{\underline{\varphi}(\underline{a})}=\mbox{dim}\,\mbox{Ker}\,\mbox{ad}^{r}T^{lk}(\underline{a})\cdot S^{lk}(\underline{a}).
\end{gather*}

\begin{proof}
The conclusion of Lemma 3.3 is the consequence of the reasoning above and Lemma 3.2.
\end{proof}

\vspace{2ex}
\textbf{Remark 3.2.} The above lemma will play a crucial role in the Section 6 for the proof that the stratification by the diagram of initial exponents implies the semicontinuity of the Hilbert-Samuel function.

\vspace{2ex}
Let $L$ be a definable leaf in $M^{s}_{\varphi}$- this means that $L$ is a connected definable subset of $M^{s}_{\varphi}$ which is Q-analytic submanifold in $M^{s}$. We define
\begin{gather*}
r^{lk}(L):=\max_{\underline{a}\in L}r^{lk}(\underline{a}),\\
\Theta^{lk}_{L}(\underline{a}):=\mbox{ad}^{r^{lk}(L)}T^{lk}(\underline{a})S^{lk}(\underline{a}),\\
d^{lk}_{L}:=\mbox{rk}\,\Theta^{lk}_{L}(\underline{a}).
\end{gather*}
Put
\begin{gather*}
d^{lk}(L):=\max_{\underline{a}}d^{lk}_{L}(\underline{a}).
\end{gather*}
Consider a set
\begin{gather*}
Y^{lk}:=\{\underline{a}\in L:\, r^{lk}(\underline{a})<r^{lk}(L)\}.
\end{gather*}
By the definition of $\Theta^{lk}_{L}(\underline{a})$, if $\underline{a}\in Y^{lk}$ then $\Theta^{lk}_{L}(\underline{a})=0$. Conversely, let $\Theta^{lk}_{L}(\underline{a})=0$. Suppose $r^{lk}(\underline{a})=r^{lk}(L)$. Then by Lemma 3.1, $(S^{lk}(\underline{a}),T^{lk}(\underline{a}))=0$, since $S^{lk}(\underline{a})\xi+T^{lk}(\underline{a})\zeta=0$ if and only if $\zeta\in\mbox{Ker}\,\mbox{ad}^{r^{lk}(\underline{a})}T^{lk}(\underline{a})S^{lk}(\underline{a})$. By the assumption, $\Theta^{lk}_{L}(\underline{a})=0$, thus $\mbox{Ker}\,\Theta^{lk}_{L}(\underline{a})$ is whole space. Therefore $\mbox{rk}\,T^{lk}(\underline{a})=0$, which leads to a contradiction. As a conclusion we obtain the fact that $\underline{a}\in Y^{lk}$ if and only if $\Theta^{lk}_{L}(\underline{a})=0$, and thus $Y^{lk}$ is a closed subset of $L$ nowhere dense in $L$. Let
\begin{gather*}
Z^{lk}:=Y^{lk}\cup\{\underline{a}\in L:\,d^{lk}(\underline{a})<d^{lk}(L)\}.
\end{gather*}
If $\underline{a}\in L\setminus Z^{lk}$, then $r^{lk}(\underline{a})=r^{lk}(L)$ and $d^{lk}(\underline{a})=d^{lk}(L)$. Let
\begin{gather*}
D^{k}:=L\setminus\bigcup_{l>k}Z^{lk}.
\end{gather*}
Then $D^{k}$ is a dense subset of $L$.

\vspace{2ex}
\textbf{Lemma 3.4(\cite{[BM-1]}).}\textit{ For all $\underline{a},\underline{a}'\in D^{k}$, $H_{\underline{a}}(k)=H_{\underline{a}'}(k)$ and  $l(\underline{a},k)=l(\underline{a}',k)$.}

\vspace{2ex}
For the proof of Lemma 3.4, see Lemma 5.3 in \cite{[BM-1]}.

\vspace{2ex}
\textbf{Definition 3.2.} \textit{The Generic Hilbert-Samuel function.} Put $H_{L}(k):=H_{\underline{a}}(k)$ and $l(L,k):=l(\underline{a},k)$, where $\underline{a}\in D^{k}$. We call $H_{L}(k)$ the generic Hilbert-Samuel function, and we call $l(L,k)$ the generic Chevalley estimate.

\vspace{2ex}
Assume that $\overline{L}$ lies in a coordinate chart for $M^{s}$. Take $\underline{a}\in L$. Let $k(\underline{a})$ be the least $k$, for which $H_{\underline{a}}(l)$ coincides with polynomial in $l$ for $l>k$. Each vertex of diagram $\n_{\underline{a}}=\n(\mathcal{R}_{\underline{a}})$ has the order less than or equal to $k(\underline{a})$. If $k\geq k(\underline{a})$ and $l\geq l(\underline{a},k)$, then $\mathcal{R}^{k}(\underline{a})=\mbox{Ker}\,\Theta^{lk}(\underline{a})$ and by Corollary 3.1

\begin{gather*}
\mbox{Ker}\,\Theta^{lk}\cap J^{k}(\underline{\varphi}(\underline{a}))^{\n_{\underline{a}}}=0,
\end{gather*}
and $\mbox{dim}\,J^{k}(\underline{\varphi}(\underline{a}))^{\n_{\underline{a}}}=\mbox{rk}\,\Theta^{lk}(\underline{a})$. Thus there exists a nonzero minor $Q(\underline{a})$ of $\Theta^{lk}(\underline{a})$ restricted to $J^{k}(\underline{\varphi}(\underline{a}))^{\n_{\underline{a}}}$ such that order of $Q(\underline{a})$ is $\mbox{rk}\,\Theta^{lk}(\underline{a})$. Let $\xi=(\xi_{\beta})_{|\beta|\leq k}$. Then by Cramer's rule
\begin{gather*}
(1)\,\,\,\Theta^{lk}(\underline{a})\xi=0
\end{gather*}
iff
\begin{gather*}
\xi_{\gamma}-\sum_{\beta\notin\Delta(\underline{a}),|\beta|\leq k}\xi_{\beta}\frac{P^{\beta}_{\gamma}(\underline{a})}{Q(\underline{a})}=0,
\end{gather*}
$\gamma\in\Delta(\underline{a}),\,|\gamma|\leq k$, where $\Delta(\underline{a}):=\mathbb{N}\setminus \n_{\underline{a}}$ and $P^{\beta}_{\gamma}(\underline{a})$ are the minors from the system of equations (1). Let $\mathcal{B}(\n_{\underline{a}}):=\{\alpha_{1},\dots,\alpha_{t}\}$ be the set of vertices of $\n_{\underline{a}}$. Then
\begin{gather*}
(2)\,\,\,G_{i}=\left(y-\underline{\varphi}(\underline{a})\right)^{\alpha_{i}}+\sum_{\gamma\in\Delta(\underline{a})}g_{i\gamma}\cdot\left(y-\underline{\varphi}(\underline{a})\right)^{\gamma},\,\,i=1,\dots,t,
\end{gather*}
is the standard basis of $\mathcal{R}_{\underline{a}}$ and
\begin{gather*}
g_{i\gamma}=\frac{P^{\alpha_{i}}_{\gamma}(\underline{a})}{Q(\underline{a})},
\end{gather*}
for all $\gamma\in\Delta(\underline{a})$, $|\gamma|\leq k$.

Let $k$ be the largest integer such that $H_{L}(l)$ coincides with a polynomial for $l\geq k$, and let $l=l(L,k)$. Thus all the vertices of $\n_{\underline{a}}$ has the order less than or equal to $k$ for $\underline{a}\in D^{k}$. Therefore the set of diagrams on $D^{k}$ is finite and thus has a minimum. Let $\underline{a}\in D^{k}$ be an element such that $\n_{\underline{a}}$ is the smallest diagram on $D^{k}$. Let
\begin{gather*}
Z_{Q}:=\{\underline{a}':\,Q(\underline{a}')=0\}.
\end{gather*}
Since $Q$ is a quasianalytic function, $Z_{Q}$ is a closed quasianalytic subset in $L$. We have the following

\vspace{2ex}
\textbf{Lemma 3.5.}\textit{ For all $\underline{a}'\in D^{k}\setminus Z_{Q}$, $\n_{\underline{a'}}=\n_{\underline{a}}$.}
\begin{proof} By (2), there exist relations $G^{i'}\in\mathcal{R}_{\underline{a}'}$, such that
\begin{gather*}
G'_{i}=\left(y-\underline{\varphi}(\underline{a}')\right)^{\alpha_{i}}+\sum_{\gamma\in\Delta(\underline{a}'),|\gamma|\leq k}g_{i\gamma}'\cdot\left(y-\underline{\varphi}(\underline{a})\right)^{\gamma}+M_{i},
\end{gather*}
where $g'_{i\gamma}=\frac{P_{\alpha}^{\alpha_{i}}(\underline{a})}{Q(\underline{a})}$ and $M_{i}\in\widehat{m}_{b'}^{k+1}$. We took $\underline{a}$ such that $\n_{\underline{a}}\leq\n_{\underline{a}'}$, thus $\alpha_{i}\in\n_{\underline{a}'}$ for each $i$. Therefore $\n_{\underline{a}}\subset\n_{\underline{a}'}$, which implicates $\n_{\underline{a}'}\leq\n_{\underline{a}}$. Finally $\n_{\underline{a}}=\n_{\underline{a}'}$, which ends the proof.
\end{proof}

\vspace{2ex}
\textbf{Definition 3.3.} We define $\n_{L}:=\n_{\underline{a}}$, for all $\underline{a}\in D^{k}\setminus Z_{Q}$. We call $\n_{L}$ the generic diagram of initial exponents for $L$.

\vspace{2ex}
\textbf{Lemma 3.6.}\textit{(\cite{[BM-1]}, Lemma 5.8)}\textit{ For all $\underline{a}\in L$, $\n_{\underline{a}}\geq \n_{L}$.}

\section{Proof of implication $(2)\Rightarrow(3)$.}

\vspace{2ex}
As we mentioned in the Introduction, to establish the semicontinuity of the diagram of initial exponents, E. Birstone and P. Milman proved simultaneously that the sets $Z_{\n}^{+}$ (see Definition 4.2 below) are subanalytic and closed, if the uniform Chevalley estimate holds (\cite{[BM-1]}, Proposition 8.6). Their proof cannot be directly applied in the quasianalytic settings, since it relies on the fact that the ring of formal power series is faithfully flat over the ring of analytic function germs, which is no longer available in the quasianalytic case. Here we are forced to follow a different, not so direct strategy. In this section, we prove that the sets $Z_{\n}^{+}$ are definable if the umniform Chevalley estimate holds (Proposition 4.3), and next the implication $(2)\Rightarrow(3)$.

We should note that Bierstone--Milman's proof of Proposition~8.6 from paper \cite{[BM-1]} contains small error. In this chapter, we give an example (Remark 4.2) which shows why their proof is not completely correct and provide a proof of Proposition 4.3 which improves their arguments.

Assume that $Z\subset X$ are closed definable sets.
Let $\mathfrak{N}\in \mathcal{D}(n)$ and $\alpha\in\mathbb{N}^{n}$. We repeat the following notation by Bierstone and Milman:

\vspace{2ex}
\textbf{Definition 4.1} Set
\begin{gather*}
\mathfrak{N}(\alpha):=\mathbb{N}^{n}+\{\beta\in\mathfrak{N}: \beta\leq\alpha\}\\
\mathfrak{N}^{-}(\alpha):=\mathbb{N}^{n}+\{\beta\in\mathfrak{N}: \beta<\alpha\}.
\end{gather*}

\vspace{2ex}
\textbf{Definition 4.2} Put
\begin{gather*}
Z_{\mathfrak{N}}(\alpha):=Z\cup\{b\in X\setminus Z:\mathfrak{N}_{b}(\alpha)\geq\mathfrak{N}(\alpha)\},\\
Z^{+}_{\mathfrak{N}}(\alpha):=Z\cup\{b\in X\setminus Z:\mathfrak{N}_{b}(\alpha)>\mathfrak{N}(\alpha)\},\\
Z_{\mathfrak{N}}:=Z\cup\{b\in X\setminus Z:\mathfrak{N}_{b}\geq\mathfrak{N}\},\\
Z^{+}_{\mathfrak{N}}:=Z\cup\{b\in X\setminus Z:\mathfrak{N}_{b}>\mathfrak{N}\}.
\end{gather*}

We have the following

\vspace{2ex}
\textbf{Lemma 4.1}\textit{ Let $\mathfrak{N}_{1}, \mathfrak{N}_{2}\in \mathcal{D}(n)$ and $\alpha\in\mathbb{N} ^{n}$. The following conditions are equivalent:}

(1) $\mathfrak{N}_{1}(\alpha)<\mathfrak{N}_{2}(\alpha)$,

(2)\textit{ there exists $\theta\leq\alpha$ such that $\mathfrak{N}_{1}^{-}(\theta)=\mathfrak{N}_{2}^{-}(\theta)$, $\theta\in\mathfrak{N}_{1}(\theta)$ and $\theta\notin\mathfrak{N}_{2}(\theta)$.}

\vspace{2ex}
Lemma 4.1. is an analogue of Remark 8.5 from \cite{[BM-1]}.

\vspace{2ex}
Now we present two crucial propositions, which we need to show that uniform Chevalley estimate implies stratification by the diagram of initial exponents. We have the following

\vspace{2ex}
\textbf{Proposition 4.1.} \textit{Assume that $X$ has the uniform Chevalley estimate relatively to $Z$ and let $Y\subset X$ be a closed definable set such that $\mathfrak{N}_{b}^{-}(\alpha)$ is constant on $Y$ for some $\alpha\in\mathbb{N}^{n}$. Then $Z\cup \{b\in Y\setminus Z : \alpha\notin\mathfrak{N}_{b}\}$ is a definable set.}

\vspace{2ex}
The proof of Proposition 4.1 is the same as the first part of the proof of Proposition 8.3 in \cite{[BM-1]}, where the authors show that the set considered is subanalytic. Since Proposition 4.1 plays important role in further reasoning, we rpeat the proof by Bierstone and Milman for readers convinience.  
\begin{proof} If $\alpha\in\mathfrak{N}_{b}^{-}(\alpha)$ then $\{b\in Y:\alpha\notin\mathfrak{N}_{b}\}=\emptyset$, therefore it is enough to consider the case, where $\alpha\notin\mathfrak{N}_{b}^{-}(\alpha)$. We assume that $X$ is compact.

Let $\varphi:M\rightarrow\mathbb{R}^{n}$ be a proper quasianalytic mapping such that $\varphi(M)=X$ (Corollary 2.1). By the assumption and Lemma 2.3, there exists a function $l_{\varphi^{*}}:\mathbb{N}\rightarrow\mathbb{N}$ such that $l_{\varphi^{*}}(b,k)\leq l_{\varphi^{*}}(k)$ for every $b\in X\setminus Z$ and $k\in\mathbb{N}$.

Let $k=|\alpha|$, $l=l_{\varphi^{*}}(k)$ and $b\in X$. Let 
\begin{gather*}
J^{l}(b)^{\mathfrak{N}^{-}(\alpha)}:=\{\xi=(\xi_{\beta})_{|\beta|\leq l}:\beta\in\mathfrak{N}^{-}(\alpha) \Rightarrow \xi_{\beta}=0 \}, 
\end{gather*}
where $\mathfrak{N}^{-}(\alpha)=\mathfrak{N}_{b}^{-}(\alpha)$ for any $b\in Y\setminus Z$. Put
\begin{gather*}
\mathfrak{N}^{-}(\alpha)_{+}:=\mathfrak{N}^{-}(\alpha)\cup \{\beta\in\mathbb{N}^{n}:\beta>\alpha\}.
\end{gather*}
We obtain a direct-sum decomposition
\begin{gather*}
J^{l}(b)^{\mathfrak{N}^{-}(\alpha)}=J^{l}(b)^{\mathfrak{N}^{-}(\alpha)_{+}}\oplus(\widehat{m}_{y}^{>\alpha}\cap J^{l}(b)^{\mathfrak{N}^{-}(\alpha)}),
\end{gather*}
where $\widehat{m}_{b}^{>\alpha}\subset \widehat{m}_{b}$ is the ideal generated by monomials $(y-b)^\beta$ for $\beta>\alpha$. For $\xi\in J^{l}(b)^{\mathfrak{N}^{-}(\alpha)}$, we write $\xi=(\eta,\zeta)$ in order to the above decomposition.

Consider $a\in\varphi^{-1}(b)$ in local coordinate chart $(x_{1},...,x_{m})$ in a neighborhood of $a$ in $M$. Then we treat $J^{l}\varphi: J^{l}(b)\rightarrow J^{l}(a)$ as a matrix and we write $(A(a),B(a))$ for the matrix of $J^{l}\varphi|_{J^{l}(b)^{\mathfrak{N}^{-}(\alpha)}}$ according to the direct-sum above. We need two lemmas(see \cite{[BM-1]}).

\vspace{2ex}
\textbf{Lemma 4.2}\textit{ Let $b\in Y\setminus Z$. Then we have}
\begin{gather*}
\alpha\notin\mathfrak{N}_{b}\Leftrightarrow \left [ \left (\forall_{a\in \varphi^{-}(b)}\ A(a)\eta+B(a)\zeta=0\right ) \Rightarrow \eta=0 \right ].
\end{gather*}

\begin{proof} Consider $\xi\in J^{l}(b)^{\mathfrak{N}^{-}(\alpha)}$. Let $P_{\xi}=\sum_{|\beta|\leq l}\xi_{\beta}(y-b)^{\beta}$ be a polynomial which generates $\xi$. If $\xi=(\eta,\zeta)$, we have $P_{\xi}=P_{\eta}+P_{\zeta}$ with respect to the direct-sum decomposition. Let $a\in\varphi^{-1}(b)$. Then $A(a)\eta+B(a)\zeta=0$ if and only if $\varphi^{*}_{a}(P_{\xi})\in m_{a}^{l+1}$.

Suppose $\alpha\notin\mathfrak{N}_{b}$ and let $\xi=(\eta,\zeta)\in J^{l}(b)^{\mathfrak{N}^{-}(\alpha)}$ and  that $A(a)\eta+B(a)\zeta=0$ for all $a\in\varphi^{-1}(b)$. Then $\varphi^{*}_{a}(P_{\xi})\in m_{a}^{l+1}$ and thus $\nu_{M,a}(\widehat{\varphi}^{*}_{a}(P_{\xi}))>l$ for all $a\in\varphi^{-1}(b)$. By Chevalley estimate and Lemma 2.3, $\nu_{X,b}(P_{\xi})>k$, therefore $P_{\xi}\in \mathcal{R}_{b}+(y-b)^{k+1}$. Since $|\alpha|\leq k$, $(y-b)^{k+1}\subset (y-b)^{>\alpha}$. According to the definition of $P_{\xi}=P_{\eta}+P_{\zeta}$ and direct-sum decomposition we see that $P_{\zeta}\in(y-b)^{>\alpha}$ and thus $P_{\eta}\in \mathcal{R}_{b}+(y-b)^{>\alpha}$. By the assumption, $\eta\in J^{l}(b)^{\mathfrak{N}^{-}(\alpha)}$. Since $\alpha \notin\mathfrak{N}_{b}$, $\mathfrak{N}_{b}(\alpha)=\mathfrak{N}^{-}(\alpha)$ and $\mathfrak{N}_{b}\cup \{\beta:\beta>\alpha\}=\mathfrak{N}^{-}(\alpha)_{+}$. Thus $P_{\eta}\in\mathbb{R}[[y-b]]^{\mathfrak{N}_{b}\cup \{\beta:\beta>\alpha\}}$. On the other hand $\mathcal{R}_{b}+(y-b)^{>\alpha}$ is an ideal in $\mathbb{R}[[y-b]]$ whose diagram is $\mathfrak{N}_{b}\cup \{\beta:\beta>\alpha\}$. By Corollary 3.1 we get

\begin{gather*}
\mathcal{R}_{b}+(y-b)^{>\alpha}\cap \mathbb{R}[[y-b]]^{\mathfrak{N}_{b}\cup \{\beta:\beta>\alpha\}}=0.
\end{gather*}

Since $P_{\eta}$ belongs to this intersection, we obtain $P_{\eta}=0$. Therefore $\eta=0$.

In order to proof inverse implication suppose tah $\alpha\in\mathfrak{N}_{b}$. By Theorem 3.1, there exists $G\in\mathcal{R}_{b}$ such that $\mbox{mon\,} G=(y-b)^{\alpha}$ and $G-(y-b)^{\alpha}\in\mathbb{R}[[y-b]]^{\mathfrak{N}_{b}}$. It is enough to show that there exists $\xi$ such that, if $(A(a),B(a))\xi=0$ for all $a\in\varphi^{-1}(b)$ then $\eta\neq 0$. Since $\mathfrak{N}^{-}_{b}(\alpha)=\mathfrak{N}^{-}(\alpha)$, $G-(y-b)^{\alpha}\in\mathbb{R}[[y-b]]^{\mathfrak{N}^{-}(\alpha)}$. Put $\xi=J^{l}G(b)$. Then $\xi\in J^{l}(b)^{\mathfrak{N}^{-}(\alpha)}$ and $(A(a),B(a))\xi=0$ for all $a\in\varphi^{-1}(b)$ since $G\in\mathcal{R}_{b}$. By the definition of $\xi=(\eta,\zeta)$, we have $\xi_{\alpha}=1$ and $\alpha\notin\mathfrak{N}^{-}(\alpha)\cup\{\beta:\beta>\alpha\}$, therefore $\eta\neq 0$.
\end{proof}

\vspace{2ex}
\textbf{Lemma 4.3}\textit{ Let $\{C(\lambda):\lambda\in\Lambda\}$ be a set of matrices from each of which has $p$ columns and $\sharp\Lambda\geq p$. Let $\mbox{Ker}\,C(\lambda)=\{\xi=(\xi_{1},\dots,\xi_{p}):C(\lambda)\xi=0\ for\ all\ \lambda\in\Lambda\}$. Then there exists $J\subset \Lambda$ such that $\sharp J=p$ and $\mbox{Ker}\,C(J)=\mbox{Ker}\,C(\lambda)$.}

\begin{proof} If $\lambda\in\Lambda$, then $C(\lambda)\xi=0$ if and only if for any row $w$ of $C(\lambda)$ the scalar product $w\cdot\xi=0$. Since the number of linearly independent rows from all $C(\lambda)$ is less than or equal to $p$, there exists such $J$.
\end{proof}
To complete the proof of proposition it is enough to show that 
\begin{gather*}
\Sigma=\{b\in Y\setminus Z:\alpha\notin\mathfrak{N}_{b}\}
\end{gather*}
is a definable set. Let $q= \left (n+l \atop l\right )$. For $\underline{a}=(a^{1},\dots,a^{q})\in M^{q}_{\varphi}$ we write

\begin{gather*}
(\underline{A}(\underline{a}),\underline{B}(\underline{a})):=\left[\begin{array}{cc} A(a^{1})&B(a^{1})\\ A(a^{2})&B(a^{2})\\ \vdots&\vdots \\ A(a^{q})&B(a^{q})  \end{array}\right ].
\end{gather*}
By Lemma 4.2, for all $b\in X$, there exits $\underline{a}\in\underline{\varphi}^{-1}(b)\subset M^{q}_{\varphi}$ such that $A(a)\eta+B(a)\zeta=0$ for all $a\in\varphi^{-1}(b)$ if and only if $\underline{A}(\underline{a})\eta+\underline{B}(\underline{a})\zeta=0$. Let $b\in Y\setminus Z$. Again, by Lemma 4.2, $b\in\Sigma$ if and only if there exists $\underline{a}\in\underline{\varphi}^{-1}(b)$ such that
\begin{gather*}
\underline{A}(\underline{a})\eta+\underline{B}(\underline{a})\zeta=0 \Longleftrightarrow \eta=0.
\end{gather*}

Let $\underline{a}\in M^{q}_{\varphi}$ and let $r(\underline{a}):=\mbox{rank}(\underline{B}(\underline{a}))$. Then $r(\underline{a})\leq q$. Put $\underline{T}=\mbox{ad}^{r(\underline{a})}\underline{B}(\underline{a})\cdot\underline{A}(\underline{a})$, then
\begin{gather*}
\mbox{Ker}\,\underline{T}(\underline{a})=\{\eta:\underline{A}(\underline{a})\eta\in\mbox{Im}\,\underline{B}(\underline{a})\}.
\end{gather*}

The set $M^{q}_{\varphi}$ can be decomposed in the following way: $M^{q}_{\varphi}=\bigcup_{r}S^{r}$, where
$S^{r}:=\{\underline{a}\in M^{q}_{\varphi}:\mbox{rank}\,\underline{B}(\underline{a})=r\}$. We observe that each $S^{r}$ is a difference of two analytic sets in $M^{q}_{\varphi}$--- $S^{r}=W^{r}\setminus W^{r-1}$, where $W^{r}:=\{\underline{a}\in M^{q}_{\varphi}: \mbox{rank}\,\underline{B}(\underline{a})\leq r\}$. For each $r$, we put $S^{r}_{0}:=\{\underline{a}\in S^{r}:\mbox{Ker}\, T=0\}$ and $S_{0}:=\bigcup_{r}S^{r}_{0}$, which is a definable set. Since $\mbox{Ker}\,T=0$ if and only if $\underline{A}(\underline{a})(\eta)+\underline{B}(\underline{a})(\zeta)=0$ implies $\eta=0$ and, by Lemma 4.2, $\Sigma=Y\setminus Z \cup \underline{\varphi}(S_{0})$, thus $\Sigma$ is definable.
\end{proof}

\vspace{2ex}
\textbf{Proposition 4.2}\textit{ Assume that $Z\subset X$ are closed definable sets such that $X$ has the uniform Chevalley estimate relatively to $Z$. Let $\mathfrak{N}\in\mathcal{D}(n)$ and $\alpha\in\mathbb{N}^{n}$. Then $Z_{\mathfrak{N}}(\alpha)$ and $Z^{+}_{\mathfrak{N}}(\alpha)$ are definable subsets of $X$. }

\vspace{2ex}
\textbf{Remark 4.1}
Proposition 4.2 is a weaker version of Proposition 8.6 from \cite{[BM-1]}. In comparison to the original proposition we cannot obtain a closedness of the set considered in the same way as in \cite{[BM-1]}, since it is unknown if $\widehat{\mathcal{Q}}_{b}$ is faithfully flat over $\mathcal{Q}_{b}$. Since there is an error in original proof we repeat with a correction which is explained in Remark 4.2. 

\vspace{2ex}
\begin{proof} We prove Proposition 4.2 by the induction on $\alpha$. First assume that $\alpha=0$. Then $\mathfrak{N}_{b}(\alpha)=\emptyset$ for $b\in X$, and there are two possible cases. If $\alpha\notin\mathfrak{N}$ then $\mathfrak{N}(\alpha)=\emptyset$, whence $Z_{\mathfrak{N}}(\alpha)=X$ and $Z^{+}_{\mathfrak{N}}(\alpha)=Z$. On the other hand, if $\alpha\in\mathfrak{N}$ then $\mathfrak{N}(\alpha)=\mathbb{N}^{n}$ and therefore $Z_{\mathfrak{N}}(\alpha)=Z^{+}_{\mathfrak{N}}(\alpha)=X$.

Now assume that the conclusion is true for all exponents $<\alpha$. Let $\beta_{b}\in\n_{b}$ be the largest element which is less then $\alpha$ for $b\in X\setminus Z$, and $\beta_{1}\in\n$ be the largest element less than $\alpha$ in $\n$. Let $\beta:=\max\{\max\{\beta_{b},b\in X\setminus Z\},\beta_{1}\}$. Consider the following sets
\begin{gather*}
X_{1}=Z_{\mathfrak{N}}(\alpha),\\
X_{0}=Z_{\mathfrak{N}}(\beta),\\
Z_{1}=Z^{+}_{\mathfrak{N}}(\alpha),\\
Z_{0}=Z^{+}_{\mathfrak{N}}(\beta).
\end{gather*}
We shall prove that $Z_{0}\subset Z_{1}\subset X_{1}\subset X_{0}$.

To show inclusion $Z_{0}\subset Z_{1}$ take $b\in\{X\setminus Z:\mathfrak{N}_{b}(\beta)>\mathfrak{N}(\beta)\}$. By Lemma 4.1, there exists $\theta\leq\beta$ such that $\mathfrak{N}^{-}_{b}(\theta)=\mathfrak{N}^{-}(\theta)$, $\theta\in\mathfrak{N}(\theta)$ and $\theta\notin\mathfrak{N}_{b}$. Since $\theta\leq\beta<\alpha$ and again by Lemma 4.1, $\mathfrak{N}_{b}(\alpha)>\mathfrak{N}(\alpha)$. Therefore 
\begin{gather*}
\{b\in X\setminus Z:\mathfrak{N}_{b}(\beta)>\mathfrak{N}(\beta)\}\subset\{b\in X\setminus Z:\mathfrak{N}_{b}(\alpha)>\mathfrak{N}_{b}(\alpha)\}
\end{gather*} 
and $Z_{0}\subset Z_{1}$.

It is clear that $Z_{1}\subset X_{1}$. To prove that $X_{1}\subset X_{0}$ it is enough to show that $\mathfrak{N}_{b}(\alpha)\geq\mathfrak{N}(\alpha)$ implies $\mathfrak{N}_{b}(\beta)\geq\mathfrak{N}(\beta)$. Suppose $\mathfrak{N}_{b}(\beta)<\mathfrak{N}(\beta)$. By Lemma 4.1, there exists $\theta\leq\beta$ such that $\mathfrak{N}^{-}_{b}(\theta)=\mathfrak{N}^{-}(\theta)$, $\theta\in\mathfrak{N}_{b}(\theta)$ and $\theta\notin\mathfrak{N}$. Since $\theta\leq\beta<\alpha$ and by Lemma 4.1, $\mathfrak{N}_{b}(\alpha)<\mathfrak{N}(\alpha)$. A contradiction, which proves that $X_{1}\subset X_{0}$.

By the induction, $X_{0}$ and $Z_{0}$ are definable sets. We prove our thesis for $\alpha$.

\textit{Case 1.} Assume $\alpha\in\mathfrak{N}^{-}(\alpha)$. Then $\n(\alpha)=\n^{-}(\alpha)=\n(\beta)$. Let $b\in X$. If $\n_{b}(\beta)=\n(\beta)$, then $\alpha\in\n_{b}(\beta)$, $\alpha\in\n(\beta)$. Also $\n_{b}(\beta)=\n_{b}(\alpha)$ and $\n(\beta)=\n(\alpha)$. Thus $\n_{b}(\alpha)=\n(\alpha)$. On the other hand, if $\n_{b}(\beta)>\n(\beta)$ then $\n_{b}(\alpha)>\n(\alpha)$. Since $\alpha\in\n_{b}\cap\n$, $\n_{b}(\alpha)\geq\n(\alpha)$ if and only if $\n_{b}(\beta)\geq\n(\beta)$, and thus $X_{1}=X_{0}$, which is definable.

If $\n_{b}(\beta)>\n(\beta)$, then by Lemma 4.1 there exists $\theta\leq\alpha$ such that $\mathfrak{N}^{-}_{b}(\theta)=\mathfrak{N}^{-}(\theta)$, $\theta\in\mathfrak{N}(\theta)$ and $\theta\notin\mathfrak{N}_{b}(\theta)$. If $\theta\leq\beta<\alpha$, then $\n_{b}(\alpha)>\n(\alpha)$. If $\beta<\theta$, then $\theta=\alpha$ by the definition of $\beta$. Thus $\n^{-}_{b}(\alpha)=\n^{-}(\alpha)$, $\alpha\in \n(\alpha)$ and $\alpha\notin\n_{b}(\alpha)$. On the other hand $\alpha\in\n^{-}(\alpha)=\n^{-}_{b}(\alpha)$, which is a contradiction. Therefore, $\n_{b}(\alpha)>\n(\alpha)$ if and only if $\n_{b}(\beta)>\n(\beta)$, and thus $Z_{1}=Z_{0}$.

\vspace{2ex}
\textit{Case 2.} Assume that $\alpha\notin\n^{-}(\alpha)$ and $\alpha\in\n$. Let $b\in X$. If $\n_{b}(\beta)=\n(\beta)$ and $\alpha\in\n_{b}$, then $\n_{b}(\alpha)\geq\n(\alpha)$, since $\alpha$ is a vertex of $\n$. If $\alpha\notin\n_{b}$ then $\n_{b}(\alpha)\geq\n(\alpha)$, since $\n(\alpha)$ has one additional vertex in comparison to $\n_{b}(\alpha)$ and $\n^{-}_{b}(\alpha)=\n^{-}(\alpha)$. Therefore $\n_{b}(\alpha)\geq\n(\alpha)$ if and only if $\n_{b}(\beta)\geq\n(\beta)$ and $X_{1}=X_{0}$.

If $\n_{b}(\alpha)>\n(\alpha)$, then there exists $\theta<\alpha$ such that $\n^{-}_{b}(\theta)=\n^{-}(\theta)$, $\theta\in\n(\theta)$ and $\theta\notin\n_{b}$. If one can find $\theta<\beta$, then $\n_{b}(\beta)>\n(\beta)$. If not, then $\theta=\alpha$ and $\n_{b}(\beta)=\n(\beta)$, $\theta\notin\n_{b}$ and $\theta\in\n$. If we assume that $\n_{b}(\beta)=\n(\beta)$, $\alpha\notin\n_{b}$ and $\alpha\in\n$, then by Lemma 4.1 $\n_{b}(\alpha)>\n(\alpha)$. Therefore $\n_{b}(\alpha)>\n(\alpha)$ if and only if either $\n_{b}(\beta)>\n(\beta)$ or $\n_{b}(\beta)=\n(\beta)$, $\alpha\in\n$ and $\alpha\notin\n_{b}$. In that case $Z_{1}=Z_{0}\cup\{b\in X_{0}\setminus Z_{0}:\alpha\notin\n_{b}\}$. $Z_{0}$ is definable by induction, $\{b\in X_{0}\setminus Z_{0}:\alpha\notin\n_{b}\}$ is definable by Proposition 4.1, and thus $Z_{1}$ is definable.

\vspace{2ex}
\textit{Case 3.} Assume $\alpha\notin\n^{-}(\alpha)$ and $\alpha\notin\n$. If $\n_{b}(\alpha)>\n(\alpha)$ then there exists $\theta\leq\alpha$ such that $\n^{-}(\theta)=\n^{-}_{b}(\theta)$, $\theta\notin\n_{b}$ and $\theta\in\n$. If $\theta\leq\beta$ then $\n_{b}(\beta)>\n(\beta)$. If not, thus $\theta=\alpha\in\n$, which is a contradiction. Therefore $\n_{b}(\alpha)>\n(\alpha)$ if and only if $\n_{b}(\beta)>\n(\beta)$ and $Z_{1}=Z_{0}$.

If $\n_{b}(\alpha)\geq\n(\alpha)$, then either $\n_{b}(\alpha)>\n(\alpha)$ or $\n_{b}(\alpha)=\n(\alpha)$. Thus $\alpha\notin\n_{b}$ and $\n_{b}(\beta)=\n(\beta)$ or $\alpha\in\n_{b}$ and $\n_{b}(\alpha)>\n(\alpha)$. Therefore $\n_{b}(\alpha)\geq\n(\alpha)$ if and only if $\n_{b}(\beta)>\n(\beta)$ or $\n_{b}(\beta)=\n(\beta)$ and $\alpha\notin\n_{b}$. Thus we obtain $X_{1}=Z_{1}\cup\{b\in X_{0}\setminus Z_{1}: \alpha\notin\n_{b}\}$, which is definable via Case 1, Case 2 and Proposition 4.1. This ends the proof.
\end{proof}

\vspace{2ex}
\textbf{Remark 4.2} In the original proof authors set $\beta$ as the largest element of $\n$, which is smaller than $\alpha$. Then, in Case 2, they claim that if $\alpha\notin\n^{-}(\alpha)$, then $\n_{b}(\alpha)>\n(\alpha)$ if and only if either $\n_{b}(\beta)>\n(\beta)$ or $\n_{b}(\beta)=\n(\beta)$ and $\alpha\in\n$, $\alpha\notin\n_{b}$. We give an example that this equivalence is not true.

Consider $\n_{b}:=\{(0,0,1),(0,1,0)\}+\mathbb{N}^{3}$, $\n:=\{(1,0,0),(0,0,1)\}+\mathbb{N}^{3}$. Let $\alpha=(1,0,0)$. It is clear that $\n_{b}(\alpha)<\n(\alpha)$. The largest element in $\n$ less then $\alpha$ is $\beta=(0,0,1)$, thus $\n_{b}(\beta)=\n(\beta)$. Of course $\alpha\in\n$, $\alpha\notin\n^{-}(\alpha)$ and $\alpha\notin\n_{b}$. Therefore the equivalence above is not valid. Now if we take $\beta$ as in our proof, then $\beta\geq(0,1,0)$, and this phenomenon does not occur.

\vspace{2ex}
\textbf{Corollary 4.1}\textit{ Let $\n\in\mathcal{D}(n)$. Then $Z^{+}_{\n}$ is a definable set.}

\begin{proof} Let $\alpha_{1}<\alpha_{2}<\dots<\alpha_{k}=\alpha$ be the vertices of $\n$ and $\beta_{1}<\beta_{2}<\dots<\beta_{l}$ be the vertices of $\n_{b}$ for $b\in X\setminus Z$. We will prove that $\n_{b}>\n$ if and only if $\n_{b}(\alpha)>\n(\alpha)=\n$.

Assume that $\n_{b}>\n$. If there exists $s\leq l$ such that $\beta_{s}>\alpha_{s}$ and $\beta_{i}=\alpha_{i}$ for $i<s$, then, since $s\leq l$, all $\beta_{i}$ and $\alpha_{i}$ for $i\leq s$ are vertices of $\n_{b}(\alpha)$ and $\n(\alpha)$ respectively. Therefore $\n_{b}(\alpha)>\n(\alpha)$. If there is no such $s$ thus for all $i\leq l$ we have $\beta_{i}=\alpha_{i}$ for $i\leq l$ and $l<k$. In that case $\n_{b}(\alpha)=\n_{b}$, and since $\n(\alpha)=\n$ we get $\n_{b}(\alpha)>\n(\alpha)$.

Now assume that $\n_{b}(\alpha)>\n(\alpha)$. If there exists $s\leq l$ such that $\beta_{s}>\alpha_{s}$ and $\beta_{i}=\alpha_{i}$ for $i<s$ we get immediately $\n_{b}>\n$. If there is no such $s$, then $\n_{b}(\alpha)=\n_{b}$ and $\n_{b}>\n$, or there exists vertex $\beta_{s}>\alpha=\alpha_{k}$, $s\leq k$ and $\n_{b}>\n$. Therefore $Z^{+}_{\n}=Z^{+}_{\n}(\alpha)$, which is definable by Proposition 4.2.

\end{proof}

\vspace{2ex}
\textbf{Theorem 4.1}\textit{ Let $Z\subset X$ be closed definable sets such that $X$ has the uniform Chevalley estimate relatively to $Z$. Then, for any compact set $K\subset X$, $\sharp\{\n_{b}:\, b\in (X\setminus Z)\cap K \}<\infty$.}

\begin{proof} We can assume that $X$ is compact definable set. Let $Y$ be a closed definable set such that $Z\subset Y\subset X$ and $\sharp\{\n_{b}:b\in(X\setminus Y)\}<\infty$. Such a set always exists because we can take as $Y$ whole $X$. We will prove that there exists closed definable set $Y'\subset Y$ such that $\mbox{dim}(Y'\setminus Z)<\mbox{dim}(Y\setminus Z)$ and $\sharp\{\n_{b}:b\in X\setminus Y'\}<\infty$.

Let $\varphi:M\rightarrow\mathbb{R}^{n}$ be a proper real analytic mapping such that $\varphi(M)=X$ and let $0\neq s\in \mathbb{N}$. We denote by $\underline{\varphi}^{s}:M_{\varphi}^{s}\rightarrow\mathbb{R}^{n}$ the induced mapping from the s-fold fibre-product and we write $\underline{x}=(x_{1},\dots,x_{s})\in M_{\varphi}^{s}$. For $\underline{x}$ we write $\varphi^{-1}(\underline{\varphi}(\underline{x}))=\bigcup_{i=1}^{r(\underline{x})}S_{i}(\underline{x})$, where $S_{i}(\underline{x})$ are the distinct connected components of $\varphi^{-1}(\underline{\varphi}(\underline{x}))$. Let $\underline{M}_{\varphi}^{s}$ be a subset of $M_{\varphi}^{s}$ such that each $x_{i}$ lies in a distinct connected component of $\varphi^{-1}(\underline{\varphi}(\underline{x}))$. Note that if $\underline{x}\in\underline{M}_{\varphi}^{s}$, then $r(\underline{x})\geq s$. By Corollary 2.7 $\underline{M}_{\varphi}^{s}$ is definable. Let $L=\varphi^{-}(Y)$ and $\underline{L}^{s}=\underline{\varphi}^{-1}(Y)\cap\underline{M}_{\varphi}^{s}$. For each $s$ we have a diagram of inclusions and projections:
\begin{gather*}
\begin{array}{ccccc}\underline{L}^{s+1}&\subset&\underline{M}^{s+1}_{\varphi}&\subset&M_{\varphi}^{s+1}\\\downarrow&&\downarrow&&\downarrow\\ \underline{L}^{s}&\subset&\underline{M}^{s}_{\varphi}&\subset&M_{\varphi}^{s}\\\downarrow&&\downarrow&&\downarrow \\ L&\subset&M&\subset&M\end{array},
\end{gather*}
where down arrows represents the projections $\pi(x_{1},\dots,x_{s+1})=(x_{1},\dots,x_{s})$.

Consider a set $\underline{L}^{s}\setminus\left((\underline{\varphi}^{s})^{-1}(Z)\cup \pi(\underline{L}^{s+1})\right)$. This is a set of those $\underline{x}\in\underline{L}^{s}$ that $\underline{\varphi}^{s}(\underline{x})\notin Z$ and $(\underline{\varphi}^{s})^{-1}(\underline{x})$ has exactly $s$ connected components. By Corollary 2.6, $\underline{L}^{s}\setminus\left((\underline{\varphi}^{s})^{-1}(Z)\cup \pi(\underline{L}^{s+1})\right)=\bigcup_{j}W_{s,j}$, where this sum is a finite partition into a smooth, connected and definable sets. Since we assumed that $X$ is compact and by Corollary 2.5, the number of connected components of the fibre $\varphi^{-1}(b)$, for $b\in X$, is bounded. Let $t$ be the largest number of connected components of the fibre. Therefore \begin{gather*}Y=Z\cup\bigcup_{s=1}^{t}\bigcup_{j}\underline{\varphi}^{s}(W_{s,j}).
\end{gather*}
For $\underline{a}\in W_{s,j}$ we have $\mathcal{R}_{\underline{a}}=\bigcap_{i=1}^{s}\mbox{Ker}\,\widehat{\varphi}^{*}_{a_{i}}=\mathcal{F}_{b}(X)$, where $\underline{\varphi}^{s}(\underline{a})=b$. Therefore
$\n_{\underline{a}}:=\n(\mathcal{R}_{\underline{a}})=\n_{b}$.
For each $s,j$, there exists a generic diagram $\n_{s,j}\in\mathcal{D}(n)$ (see Chapter 2) such that $\n_{\underline{a}}=\n_{s,j}$ on an open and dense set in $W_{s,j}$, and $\n_{\underline{a}}\geq\n_{s,j}$. By Corollary 4.1, $Z^{+}_{\n_{s,j}}$ is definable set, thus $\underline{\varphi}^{s}(W_{s,j})\setminus Z^{+}_{\n_{s,j}}$ is definable. For all $b\in\underline{\varphi}^{s}(W_{s,j})\setminus Z^{+}_{\n_{s,j}}$ we have $\n_{b}=\n_{s,j}$ and for all $b\in \underline{\varphi}^{s}(W_{s,j})\cap Z^{+}_{\n_{s,j}}$ we have $\n_{b}>\n_{s,j}$. We will prove that $\underline{\varphi}^{s}(W_{s,j})\setminus Z^{+}_{\n_{s,j}}$ is dense in $\overline{\underline{\varphi}^{s}(W_{s,j})}$.

Let $A_{s,j}$ be a subset of $W_{s,j}$ such that for all $\underline{a}\in A_{s,j}$, $\n_{\underline{a}}>\n_{s,j}$. Then $\underline{\varphi}^{s}(W_{s,j}\setminus A_{s,j})=\underline{\varphi}^{s}(W_{s,j})\setminus Z^{+}_{\n_{s,j}}.$ Since $\underline{\varphi}^{s}$ is continuous and proper map we have $\underline{\varphi}^{s}(\overline{W_{s,j}\setminus A_{s,j}})=\underline{\varphi}^{s}(\overline{W_{s,j}})=\overline{\underline{\varphi}^{s}(W_{s,j})\setminus Z^{+}_{\n_{s,j}}}$. Now let
\begin{gather*}
Z_{s,j}:=\overline{\underline{\varphi}^{s}(W_{s,j})}\setminus\left(\underline{\varphi}^{s}(W_{s,j})\setminus Z^{+}_{\n_{s,j}}\right).
\end{gather*}
Consider $Y'=Z\cup \bigcup_{s,j}\overline{Z_{s,j}}$. Since $\mbox{dim}\,Z_{s,j}<\mbox{dim}\,Y$ for each $s,j$ we have $\mbox{dim}\,Y'\setminus Z<\mbox{dim}\,Y\setminus Z$. Finally $Z\subset Y'\subset Y$ and $Y'$ is definable.
\end{proof}

\vspace{2ex}
\textbf{Corollary 4.2}\textit{ If $X$ has the uniform Chevalley estimate relatively to $Z$, then $Z_{\n}$ is definable set. }
\begin{proof}
We can adapt here the proof of Corollary 8.9 from \cite{[BM-1]}. We can assume that $X$ is compact. Let $\alpha$ be greater or equal to the largest vertex of $\n$ and $\n_{b}$ of all $\n_{b}$, $b\in X\setminus Z$. By Theorem 4.1, there exists such $\alpha$. Therefore $\n_{b}\geq\n$ if and only if $\n_{b}(\alpha)\geq\n(\alpha)$, and thus $Z_{\n}=Z_{\n}(\alpha)$.
\end{proof}

\vspace{2ex}
\textbf{Corollary 4.3}\textit{ Let $Z\subset X$ be closed definable sets, such that $X$ has the uniform Chevalley estimate relatively to $Z$. Then there exists a stratification of $X$ such that the diagram of initial exponents is constat on each stratum and $Z$ is a sum of strata.}
\begin{proof}
We can assume that $X$ is compact. Then there is a finite number of diagrams $\n_{1},\dots,\n_{k}$ of initial exponents on $X\setminus Z$. Since $Z_{\n_{i}}$ and $Z_{\n_{i}}^{+}$ are definable for $i=1,\dots,k$, thus $Y_{\n_{i}}:=Z_{\n_{i}}\setminus Z_{\n_{i}}^{+}$ is also definable. On the other hand $Y_{\n_{i}}=\{b\in X\setminus Z:\,\n_{b}=\n_{i}\}$. Thus we can write
\begin{gather*}
X=Z\cup\bigcup_{i=1}^{k}Y_{\n_{i}}.
\end{gather*}
Since there is a stratification of $X$ compatible with $Z$ and $\{Y_{\n_{i}}\}_{i=1,\dots,k}$, and the diagram of initial exponents is constant on each $Y_{\n_{i}}$ we obtain desired stratification.
\end{proof}

We proved that the uniform Chevalley estimate implies a stratification by the diagram of initial exponents. We shall prove that reverse implication also holds. We have

\vspace{2ex}
\textbf{Theorem 4.2}\textit{ Suppose that $X$ admits a stratification by the diagram of initial exponents such that $Z$ is the union of strata. Then there $X$ has an uniform Chevalley estimate to $Z$.}

\vspace{2ex}
The proof of Theorem 4.2 is the consequence of the two propositions below, which are quasianalytic analogues of Proposition 8.14 and proposition 8.15 from \cite{[BM-1]}. Let $X$ be a closed definable subset of $\mathbb{R}^{n}$.

\vspace{2ex}
\textbf{Proposition 4.3 (\cite{[BM-1]}, Proposition 8.14)}\textit{ Suppose that $Y\subset X$ is a definable set such that $\n_{b}$ is constant on $Y$. Let $K\subset X$ be compact. Then there exists $l_{K}(k)$ such that $l_{\varphi_{*}}(b,k)\leq l_{K}(k)$ for all $b\in K\cap Y$.}

\vspace{2ex}
The proof of Proposition 4.3 can be easily reduced to the following

\vspace{2ex}
\textbf{Proposition 4.4 (\cite{[BM-1]}, Proposition 8.15)}\textit{ Let $s\geq1$ and $\underline{\varphi}:M_{\varphi}^{s}\rightarrow\mathbb{R}^{n}$. Let $L$ be a relatively compact definable subset of $M_{\varphi}^{s}$ such that $\n_{\underline{a}}=\n(\mathcal{R}_{\underline{a}})$ is constant on $L$. Then there exists $l_{L}(k)$ such that $l_{\varphi^{*}}(\underline{a},k)\leq l_{L}(k)$, for $\underline{a}\in L$.}

\vspace{2ex}
The proof of Proposition 4.4 is based on the analysis of jets and the systems of linear equations which coefficients are Q-analytic functions. The reason why it can be carried over from the analytic case is the fact that the analysis mentioned above are reduced to the properties of the ring of formal power series and several good properties which Q-analytic functions share with analytic functions, for instance the property of identity. Therefore we could just repeat the proof by E. Bierstone and P. Milman.

\section{Proof of implication $(2)\Rightarrow(4)$.}

\vspace{2ex}
In this chapter we shall prove that the uniform Chevalley estimate implies the Zariski semi-continuity of the diagram of initial exponents. Let $Z\subset X$ be closed definable sets in $\mathbb{R}^{n}$ such that $X$ has the uniform Chevalley estimate relatively to $Z$. Let $\varphi:M\rightarrow\mathbb{R}^{n}$ be a proper Q-analytic mapping from a Q-analytic manifold $M$ such that $\varphi(M)=X$. By Theorem 4.1, for any compact $K$, the set of diagrams of initial exponents is finite on $K\cap X$. Therefore, to prove Zariski semi-continuity it is sufficient to prove that the set
\begin{gather*}
Z_{\n}:=Z\cup\{b\in X\setminus Z:\,\n_{b}\geq\n\}
\end{gather*}
is a closed definable set. By Corollary 4.2, $Z_{\n}$ is definable, thus it remains to prove that $Z_{\n}$ is closed.

In the paper \cite{[BM-1]} the authors proved that in the classical analytic case $Z_{\n}$ is a closed subanalytic set. They used faithful flatness of the ring of formal power series over the ring of germs of analytic functions to deduce an existance of convergent solution to system of linear equations with analytic coefficients under assumtion that there is a formal solution. In our reasoning we are forced to provide a different method, since it is not known if the ring $\mathbb{R}[[y-b]]$ of formal power series is faithfully flat over the ring $\mathcal{Q}_{b}$ of germs of Q-analytic functions at $b$. This is an open problem related to the problem of noetherianity of $\mathcal{Q}_{b}$, which has been widely studied for the past several decades, but remains unsolved.

Our proof relies on a reduction to the analysis only of definable arcs, and on Proposition 5.1, which is a special case of Proposition 8.3 from \cite{[BM-1]} for closed definable arcs.

In order to prove Proposition 5.1, we introduce a concept of an essential point $\underline{a}\in M_{\varphi}^{q}$ which determines the diagram $\n_{b}$ with $\underline{\varphi}(\underline{a})=b$ (Definition 5.2). We show that the set of essential points is definable. Then we apply curve selection to find a definable arc $\tilde{L}$ lying over $L$ and contained in the set of essential points. Next, by Puiseux's theorem, we are able to reduce the proof to the analysis of jets parameterized by Q-analytic functions of one variable. Since the local rings of quasianalytic functions of one variable are noetherian, we can find a quasianalytic solution to a system of linear equations, which describes when a multi-index $\alpha$ belongs to the diagram of initial exponents (Proposition 5.1).

\vspace{2ex}
\textbf{Definition 5.1.} Let $M$ be a Q-analytic manifold. We say that $l:[0,\epsilon]\rightarrow M$ is a definable arc if $l$ is a continuous, definable and injective function.

\vspace{2ex}
\textbf{Remark 5.1.} It follows from the cell decomposition that a closed definable subset of pure dimension 1 is a finite sum of images of definable arcs.
By abuse of terminology, by a definable arc we often mean both the arc and its image $L=l([0,\epsilon])$. In order to prove Proposition 5.1, we need a quasianalytic version of Puiseux's theorem stated below. It is a special case of the quasianalytic version of Puiseux's theorem with parameter due to K.J. Nowak(\cite{[KN-3]}).

Let us notice that subanalytic arcs are analytic curves and their local analytic rings are noetherian. Yet the former is no longer true in quasianalytic structures, as shown by K.J. Nowak in the example constructed in paper \cite{[KN-12]}. The latter seems to be doubtful as well, being related to the failure of the following splitting problem posed by K.J. Nowak in
papers \cite{[KN-7]} and \cite{[KN-5]}:

\vspace{2ex}
\textit{Let $f$ be a Q-analytic function at $0 \in \matR^{k}$ with Taylor series $\widehat{f}$. Split the set $\matN^{k}$ of exponents into two
disjoint subsets $A$ and $B$, $\matN^{k} = A \cup B$, and decompose the formal series $\widehat{f}$ into the sum of two formal
series $G$ and $H$, supported by $A$ and $B$, respectively. Do there exist two Q-functions $g$ and $h$ at $0\in \matR^{k}$ with
Taylor series $G$ and $H$, respectively?}

\vspace{2ex}
In some special cases of splitting the Taylor exponents, a negative
answer was given by H. Sfouli \cite{[Sf]}.

Although definable arc not need to be Q-analytic curve, it can be parameterized by a Q-analytic function. It is a consequence of the following

\vspace{2ex}
\textbf{Theorem 5.1}\textit{(Puiseux's Theorem). Let}
\begin{gather*}
f:(0,1)\rightarrow\mathbb{R}
\end{gather*}
\textit{be a bounded definable function. Then there exists an interval $I:=(-\epsilon,\epsilon)$ such that}

\textit{(1) either the function $f$ vanishes on $I\cap(0,1)$;}

\textit{(2) or there exist $r\in\mathbb{N}$ and a definable function $F(t)$, Q-analytic on $I$ such that}
\begin{gather*}
f(t)=F(t^{1/r}),\,\,\,for\,all\,t\in I\cap(0,1).
\end{gather*}

Our goal is the following

\vspace{2ex}
\textbf{Theorem 5.2}\textit{ The set $Z_{\n}$ is a closed definable set. }
\begin{proof}
We can assume that $X$ is compact. By Theorem 4.1 the set of the diagrams on $X\setminus Z$ is finite. By Corollary 4.3, $Z_{\n}$ is a definable set. Thus it remains to prove that $Z_{\n}$ is closed.

Let $\mathcal{W}$ be a stratification of $X$ such that $Z_{\n}$ is a sum of strata and the diagram of initial exponents is constant on each stratum. Let $W\in\mathcal{W}$. It is enough to show, that
\begin{gather*}
(*)\,\,\,\n_{b}\geq\n_{W}\,\,\,\text{for each}\,\,b\in W'
\end{gather*}
where $\n_{W}$ is the diagram of initial exponents of the ideal $\mathcal{F}_{y}(X)$ for all $y\in W$, and $W'\subset \overline{W}\setminus W$ is a stratum such that $W'\cap Z=\emptyset$. Take $b\in W'$. By the curve selection lemma (\cite{[LvDr]},Chap. 6, Corollary 1.5) there exists a continuous definable, injective function $f:(0,\epsilon)\rightarrow W$ such that $\lim_{x\rightarrow 0}f(x)=b$. It is then enough to prove $(*)$ for the definable arcs with the end at $W'$. Theorem 5.2 will be proved once we establish the following

\vspace{2ex}
\textbf{Proposition 5.1.}\textit{ Let $L\subset X$ be a closed definable arc. Let $\alpha\in\mathbb{N}^{n}$, and for each $b\in L$, $\n_{b}(\alpha)^{-}=\n(\alpha)$. Then the set $Z\cup\{b\in L\setminus Z:\,\alpha\notin\n_{b}\}$ is closed definable set.}

\vspace{2ex}
Indeed, let us assume, that for some definable arc $L$ with end $b\in\ W'$ we have $\n_{b}<\n_{W}$. Let $\beta_{1},\dots,\beta_{s}$ be the vertices of $\n_{b}$, and let $\gamma_{1},\dots,\gamma_{r}$ be the vertices of $\n_{W}$. Clearly $\n_{b}<\n_{W}$ if and only if $r<s$ and, for $i\leq r$, $\beta_{i}=\gamma_{i}$, or there exists $j\leq\min\{s,r\}$ such that $\beta_{j}<\gamma_{j}$ and $\beta_{i}=\gamma_{i}$ for $i<j$. Therefore, there exists $\alpha\in\n_{b}$ such that $\n_{b}(\alpha)^{-}=\n_{W}(\alpha)^{-}$. Then $Z\cup\{b'\in L:\alpha\notin\n_{b'}\}=Z\cup L\setminus\{b\}$, which is not closed. This contradicts the conclusion of Proposition 5.1.

\begin{proof}(of Proposition 5.1.)
By Proposition 4.2, $Z\cup\{b\in L:\alpha\notin\n_{b}\}$ is definable. It remains to prove, that it is a closed set. Let $\Sigma$ be the set from the proof of Proposition 4.1 for $Y=L$. To prove that $Z\cup\Sigma$ is closed it is enough to show, that $(L\setminus Z)\setminus\Sigma$ is open in $L\setminus Z$. If $\alpha\in\n$, then $(L\setminus Z)\setminus\Sigma$ is an empty set and thus open. Whence we can assume that $\alpha\notin\n$ and, since we only will consider multi-indices smaller then or equal to $\alpha$, we can assume that $\n=\n^{-}(\alpha)$. Therefore  $\n\subset\n_{b}$ for all $b\in L\setminus Z$.

Let $k=|\alpha|$ and $l=l_{\varphi^{*}}(k)$ (where $l_{\varphi^{*}}(k)$ is the uniform Chevalley estimate). Let $q$ be the largest number of connected components of $\varphi^{-1}(b)$ for $b\in X$. Fix a point $b\in (L\setminus Z)\setminus\Sigma$. Obviously we can assume that $b$ is the end of $L$. Let $\tau:[0,1]\rightarrow L$ be the parametrization of $L$ such that $\tau(0)=b$.

\vspace{2ex}
In order to prove Proposition 5.1, we introduce the concept of essential point.

\vspace{2ex}
\textbf{Definition 5.2.} A point $\underline{a}\in M_{\varphi}^{q}$ is called the essential point if the following implication is true: if $\xi=(\eta,\zeta)\in\mbox{Ker}\,\Phi(\underline{a})$, then $\xi\in\mbox{Ker}\,(A(a'),B(a'))$ for all $a'\in\varphi^{-1}(\underline{\varphi}(\underline{a}))$ (here $A(a')$ and $B(a')$ are matrices as in Chapter 4).

\vspace{2ex}
We have the following

\vspace{2ex}
\textbf{Lemma 5.1.}\textit{ The set $S\subset M_{\varphi}^q$ of essential points is a definable set and $\underline{\varphi}(S)=X\setminus Z$.}

\begin{proof}
Clearly $S=\{\underline{a}\in M_{\varphi}^{q}:\, \forall_{\underline{a}'}\in M_{\varphi}^{q}\, \mbox{rank}\,\Phi(\underline{a})\geq \mbox{rank}\,\Phi(\underline{a}')\}$ thus it is definable set. Finally, $\underline{\varphi}(S)=X\setminus Z$ by Lemma 4.3 for $p=q={n+l\choose l}$. 
\end{proof}

Consider a sequence $\{b_{\omega}\}_{\omega\in\mathbb{N}}\subset L$ with $\lim_{\omega\rightarrow\infty}b_{\omega}=b$. Such a sequence is a relatively compact set, because, by the assumption on $X$, $L$ is compact. Since $\underline{\varphi}$ is a proper map, the set $\underline{\varphi}^{-1}(\{b_{\omega}\}_{\omega\in\mathbb{N}}\cup\{b\})$ is a compact subset of $M_{\varphi}^{q}$. By Lemma 5.1, $\underline{\varphi}(S)=X\setminus Z$, and thus we can take a sequence $\{\underline{a}_{n}\}_{n\in\mathbb{N}}\subset S$ such that $\underline{\varphi}(\underline{a}_{\omega})=b_{\omega}$. Of course $\{\underline{a}_{\omega}\}_{\omega\in\mathbb{N}}$ has an accumulation point at $\underline{a}\in\overline{S}$ such that $\underline{\varphi}(\underline{a})=b$. Therefore, there exists a definable arc $\tilde{L}=\tilde{\tau}([0,1))$, where $\tilde{\tau}:[0,1)\rightarrow M_{\varphi}^{q}$ is a parametrization of $\tilde{L}$, such that $\tilde{\tau}(0)=\underline{a}$ and $\tilde{\tau}((0,1))\subset S$. Consider the following diagram
\begin{gather*}
\xymatrix{
\underline{\varphi}\circ \tilde{\tau}:[0,1) \ar[rd]_{\tau^{-1}\circ\underline{\varphi}\circ \tilde{\tau}} \ar[r]\,\,&
L \\&[0,1)\ar[u]_\tau}
\end{gather*}
By the monotonicity theorem (\cite{[LvDr]},Chap. 3, Theorem 1.2) and the fact that $\tau^{-1}\circ\underline{\varphi}\circ \tilde{\tau}$ is not constant, we can assume that $\tau^{-1}\circ\underline{\varphi}\circ \tilde{\tau}$ is a strictly increasing function on some interval $[0,p)\subset[0,1)$. Therefore $\underline{\varphi}(L')\supset \tau([0,p))$. By Puiseux's theorem, there exists a Q-analytic parametrization \begin{gather*}
\varepsilon:(-1,1)\rightarrow M_{\varphi}^{q}
\end{gather*}
such that $\varepsilon([0,1))=\tilde{\tau}([0,p))$.

Since $b\in(L\setminus Z)\setminus\Sigma$, $\alpha\in\n_{b}$. Thus there exists $G\in\mathcal{R}_{b}$ such that $\mbox{mon}\,G=(y-b)^{\alpha}$ and $G-(y-b)^{\alpha}\in\widehat{\mathcal{Q}}_{b}^{\n_{b}}\subset\widehat{\mathcal{Q}}_{b}^{\n}$. Here we identify $\widehat{\mathcal{Q}}_{b}$ with $\mathbb{R}[[y-b]]$. Since $\alpha\notin\n^{-}(\alpha)=\n$, we get $G\in\widehat{\mathcal{Q}}_{b}^{\n}$, and therefore $D^{\beta}G=0$ for $\beta\in\n$. Of course $J^{l}_{b}G\in J^{l}(b)\otimes_{\mathbb{R}}\widehat{\mathcal{Q}}_{b}$.

Let $\underline{a}=\varepsilon(0)$. The mapping $\underline{\varphi}:M_{\varphi}^{q}\rightarrow\mathbb{R}^{n}$ induces a ring homomorphism $\widehat{\underline{\varphi}}_{\underline{a}}^{*}:\widehat{\mathcal{Q}}_{b}\rightarrow\widehat{\mathcal{Q}}_{M_{\varphi}^{q},\underline{a}}$, and further the homomorphism
\begin{gather*}
 (1)\,\,\,\,\,J^{l}_{\underline{a}}\varphi:J^{l}(b)\otimes\widehat{O}_{M^{q}_{\varphi},\underline{a}}\rightarrow\bigoplus_{i=1}^{q}J^{l}(a_{i})\otimes\widehat{O}_{M^{q}_{\varphi},\underline{a}}.
\end{gather*}
Let $\widehat{\xi}_{\underline{a}}:=(J^{l}_{b}G)\circ\widehat{\underline{\varphi}}_{\underline{a}}$. Thus $\widehat{\xi}_{\underline{a}}\in J^{l}(b)^{\n}\otimes\widehat{O}_{M_{\varphi}^{q},\underline{a}}$ and we can write $\widehat{\xi}_{\underline{a}}=(\widehat{\eta}_{\underline{a}},\widehat{\zeta}_{\underline{a}})$ according to the direct decomposition
\begin{gather*} J^{l}(b)^{\mathfrak{N}(\alpha)}=J^{l}(b)^{\mathfrak{N}^{-}(\alpha)_{+}}\oplus(\widehat{m_{y}}^{>\alpha}\cap J^{l}(b)^{\mathfrak{N}(\alpha)}).
 \end{gather*}
Therefore the $\alpha^{th}$ component of $\widehat{\eta}_{\underline{a}}$ is 1 (since $D^{\alpha}G=1$).

The restriction of $J^{l}_{\underline{a}}\varphi$ to $J^{l}(b)^{\n}\otimes\widehat{O}_{M^{q}_{\varphi},\underline{a}}$ can be interpreted as a matrix $\Phi_{\underline{a}}=(\underline{A},\underline{B})$
with entries in $\mathcal{Q}_{M^{q}_{\varphi},\underline{a}}\subset\widehat{\mathcal{Q}}_{M^{q}_{\varphi},\underline{a}}$. Thus
we get $\Phi(\underline{a})=(\underline{A}(\underline{a}),\underline{B}(\underline{a}))$.

\vspace{2ex}
\textbf{Observation.}
If $G\in\mathcal{R}_{\underline{a}}$, we have
\begin{gather*}
J_{\underline{a}}^{l}\varphi(J_{b}^{l}G)=0.
\end{gather*}
It is an immediate consequence of the following formula
\begin{gather*}
(2)\,\,\,J^{l}_{a}\varphi\left((\widehat{\varphi}_{a}^{*}(D^{\beta}G))_{|\beta|\leq l}\right)=\left(D^{\alpha}\left(\widehat{\varphi}^{*}_{a}(G)\right)\right ) _{|\alpha|\leq l},
\end{gather*}
for all $a_{i}$ such that $\underline{a}=(a_{1},\dots,a_{q})$.

\vspace{2ex}
Hence $\widehat{\xi}_{\underline{a}}$ is a solution to the system of linear equations with Q-analytic coefficients at $\underline{a}\in M^{q}_{\varphi}$ that corresponds to the matrix $\Phi_{\underline{a}}$:
\begin{gather*}
(\star)\,\,\,\Phi_{\underline{a}}\cdot\widehat{\xi}_{\underline{a}}=0.
\end{gather*}

We shall consider the numerical system of linear equations obtained form $(\star)$ by evaluating its coefficients at point $\underline{a}'$ of arc $\tilde{L}$ near $\underline{a}$. Our goal is to find solutions $\xi(\underline{a}')$, $\underline{a}'\in\tilde{L}$, whose $\alpha$-th component $\xi^{\alpha}(\underline{a}')\neq0$.

To this end, consider the pull-back of system $(\star)$ by mean of the parametrization $\varepsilon(t)$ of the arc $\tilde{L}$:
\begin{gather*}
\Phi_{0}\cdot\widehat{\xi}_{0}=0,\,\,\text{with}\,\,\Phi_{0}=\widehat{\varepsilon}^{*}_{0}(\Phi_{\underline{a}}), \,\,\widehat{\xi}_{0}=\widehat{\varepsilon}^{*}_{0}(\widehat{\xi}_{\underline{a}})\in J^{l}(b)\otimes\widehat{Q}_{1}.
\end{gather*}
The coefficients of the system of linear equations obtained belong to the quasianalytic local ring $(Q_{1},m)$, which is a discrete valuation ring, and therefore a noetherian ring. Hence and by Lemma 2.6, there exists a solution $\xi\in J^{l}(b)^{\n}\otimes Q_{1}$ such that $\widehat{\xi}-\xi\in J^{l}(b)\otimes\widehat{m}$. Therefore, since $\widehat{\xi}^{\alpha}(0)=1$, we get $\xi^{\alpha}(0)=1$ and $\xi^{\alpha}(t)\neq 0$ for $t$ close to $0$.

In this manner, we achived numerical solutions $\xi(\varepsilon(t)):=\xi(t)$ of the system $(\star)$ at points $\underline{a}'=\varepsilon(t)$ lying on the arc $\tilde{L}$ near $\underline{a}$. At this stage, we are going to complete the proof.

Take a polynomial $f$ on $\mathbb{R}^{n}$ such that $J^{l}f(b')=\xi(\underline{a}')$ where $\underline{\varphi}(\underline{a}')=b'$. Then $\varphi^{*}_{\underline{a}'}(f)\in m^{l+1}_{\underline{a}'}$ for all $a'\in\varphi^{-1}(b')$. Hence, by the uniform Chevalley estimate 
\begin{gather*}
l_{\varphi^{*}}(b',k)\leq l_{\varphi^{*}}(k)=l,
\end{gather*}
there exists $g\in\mathcal{R}_{b'}$ such that $f-g\in\widehat{m}_{b'}^{k+1}$. Then $J^{l}g(b')=(\eta(\underline{a}'),\zeta)$ with some component $\zeta$, and $\eta^{\alpha}(\underline{a}')\neq 0$. Since $g\in\mathcal{R}_{b'}$ and $\n_{b'}^{-}(\alpha)=\n$, $\eta^{\alpha}(\underline{a}')$ is the only nonzero component of $\eta(\underline{a}')$. Therefore, $\mbox{exp}\,g=\alpha$ and $\alpha\in\n_{b'}$. Consequently, $(L\setminus Z)\setminus\Sigma$ is open in $L\setminus Z$, whence $Z\cup\Sigma$ is a closed subset, as asserted in Proposition 5.1
\end{proof}
This completes the proof Theorem 5.2.
\end{proof}

\section{The Hilbert-Samuel function.}

Here we provide a proof of the implication $(3)\Rightarrow(5)$, which is
based on  ideas similar to those from our proof of the implication $(2)\Rightarrow(4)$. Let us emphasize that Bierstone--Milman's proof from \cite{[BM-1]} does not work in the quasianalytic settings, since they use the fact that subanalytic arcs are analytic curves and their local analytic rings are noetherian. As we mentioned in the previous section, it is not true in the quasianalytic settings.

Let $Z\subset X$ be closed definable subsets of $\mathbb{R}^{n}$. For $b\in X$, $\n_{b}$ is the diagram of initial exponents of the ideal $\mathcal{R}_{b}=\mathcal{F}_{b}(X)$, and $H_{b}$ denote the Hilbert-Samuel function of $\widehat{\mathcal{Q}}_{b}/\mathcal{R}_{b}$:
\begin{gather*}
H_{b}(k)=\mbox{dim}_{\mathbb{R}}\,\widehat{\mathcal{Q}}_{b}/(\mathcal{R}_{b}+\widehat{m}_{b}^{k+1}),\,\,k\in\mathbb{N}.
\end{gather*}
The set of Hilbert-Samuel functions is equiped with the standard partial ordering, i.e. for two such functions H and H', $H \leq H'$ if $H(k) \leq H'(k)$ for all $k \in N$. With respect to the ordering above, we have the following

\vspace{2ex}
\textbf{Theorem 6.1}\textit{ Assume that $X$ admits a definable stratification such that $Z$ is the sum of strata an the diagram $\n_{b}$ is constant on each stratum disjoint with $Z$. Then $H_{b}$ is Zariski-semicontinuous relatively to $Z$.}

\begin{proof} As in the proof of Theorem 4.2, we can assume that $X$ is compact. Let $b\in X$ and let $\n_{b}$ be the diagram of initial exponents of $\mathcal{F}_{b}(X)$. By Corollary 3.2, $H_{b}(k)=\sharp\{\gamma\in\mathbb{N}^{n}\setminus\n_{b}:\,|\gamma|\leq k\}$. It follows from the stratification by the diagram of initial exponents that the function $b\rightarrow H_{b}$ is constant on each stratum. Let $W$ be a stratum disjoint with $Z$ such that the diagram of initial exponents is constant on $W$ and let $H_{W}$ be the Hilbert-Samuel function on $W$. It is sufficient to prove that for each definable arc $L$, whose interior is contained in $W$ and the end $b$ of $L$ belongs to $(\overline{W}\setminus W)$, $H_{b}(k)\geq H_{W}(k)$ for each $k\in\mathbb{N}$.

Let $q$ be the maximal number of connected components of $\varphi^{-1}(b)$ for $b\in X$. Let us consider the mapping $\underline{\varphi}:M_{\varphi}^{q}\rightarrow X$. By Corollary 2.4, there is a stratification of $M_{\varphi}^{q}=\bigcup_{i} M_{i}^{q}$ such that $\n_{\underline{a}}=\n(\mathcal{R}_{\underline{a}})$ is constant on each $M_{i}^{q}$, and this stratification is compatible with stratification by the diagram of initial exponents. For each $\underline{a}\in M_{\varphi}^{s}$ such that $\underline{\varphi}(\underline{a})=b$ there is $\mathcal{R}_{b}\subset\mathcal{R}_{\underline{a}}$, and thus $H_{b}(k)\geq H_{\underline{a}}(k)$. By Proposition 4.4, there is a uniform Chevalley estimate $l(k)$ on a stratum $M_{i}^{q}$ such that $l_{\varphi^{*}}(\underline{a},k)\leq l(k)$. Whence it is sufficient to prove the following

\vspace{2ex}
\textbf{Lemma 6.1}\textit{ Let $\tilde{L}$ be a definable arc in $M_{\varphi}^{q}$ and let $\underline{a}$ be the one of its ends. Suppose that $H_{\underline{x}}$ is constant on $\tilde{L}\setminus\{\underline{a}\}$ and $\l_{\varphi^{*}}(\underline{a},k)\leq l(k)$. Then $H_{\underline{x}}(k)\leq H_{\underline{a}}(k)$.}
\begin{proof}
Let $\varepsilon:(-1,1)\rightarrow M^{q}$ be a quasianalytic parametrization of $\tilde{L}$. Let $k\in\mathbb{N}$ and $l=l(k)$. Let $v_{1},\dots,v_{s}$ be a basis of the vector space
\begin{gather*}
 \left(\mathcal{R}_{\underline{a}}+\widehat{m}_{\underline{\varphi}(\underline{a})}^{k+1}\right)/\widehat{m}_{\underline{\varphi}(\underline{a})}^{k+1}.
\end{gather*}
Let $G_{j}\in\mathcal{R}_{\underline{a}}$ be a representation of $v_{1},\dots,v_{s}$. By $\xi_{j}^{l}$ we denote the elements of $J^{l}(\underline{\varphi}(\underline{a}))\otimes \widehat{O}_{M_{\varphi}^{q},\underline{a}}$ such that $\xi_{j}=(J^{l}_{b}G)\circ\widehat{\underline{\varphi}}_{\underline{a}}$ for $j=1,\dots,s$. By the notation from  Chapter 3, $\xi_{j}^{l}=(\xi_{j}^{k},\xi_{j}^{lk})$. For each $\underline{x}\in\tilde{L}$, we write $J^{l}\varphi(\underline{x})$ as a block matrix $J^{l}\varphi(\underline{x})=(S^{lk}(\underline{x}),T^{lk}(\underline{x}))$. Clearly, we can assume that the rank of $T^{lk}(\underline{x})$ is constant on $\tilde{L}\setminus\{\underline{a}\}$. Put $r:=\mbox{rank}\,T^{lk}(\underline{x})$. By the Observation from Chapter 5, we have
\begin{gather*}
J_{\underline{a}}^{l}\varphi\cdot\xi_{j}^{l}=0,\,\,\,j=1,\dots,s,
\end{gather*}
where $J_{\underline{a}}^{l}\varphi=(S_{\underline{a}}^{lk},T_{\underline{a}}^{lk})$ is a matrix with coefficients from $\widehat{\mathcal{Q}}_{M^{q}_{\varphi},\underline{a}}$. Then
\begin{gather*}
(\star)\,\,\,\,\,\mbox{ad}^{r}T_{\underline{a}}^{lk}\cdot S_{\underline{a}}^{lk}\cdot\xi_{j}^{k}=0,\,\,\,j=1,\dots,s.
\end{gather*}

Similarly to the proof of Proposition 4.1, we consider the pull-back of the system $(\star)$ by mean of parametrization $\varepsilon(t)$, and thus we obtain the system of linear equations with coefficients from $Q_{1}$:
\begin{gather*}
\mbox{ad}^{r}T_{0}^{lk}\cdot S_{0}^{lk}\cdot\rho_{j}^{k}=0,\,\,\,j=1,\dots,s,
\end{gather*}
where $T_{0}=\widehat{\varepsilon}^{*}_{0}(T_{\underline{a}})$, $S_{0}=\widehat{\varepsilon}^{*}_{0}(S_{\underline{a}})$ and $\rho_{j}^{k}=T_{0}=\widehat{\varepsilon}^{*}_{0}(\xi_{j}^{k})$ for $j=1,\dots,s$.

Let $w_{1},\dots,w_{p}$ be the system of generators of $\mbox{Ker}\,\mbox{ad}^{r}T_{0}^{lk}\cdot S_{0}^{lk}$. Since $\xi_{j}^{k}(0)=v_{j}$ are linearly independent, $w_{1},\dots,w_{p}$ span a vector space of dimension $\geq s$. Since $Q_{1}$ is noetherian, there exists a system of generators $w_{1}(t),\dots,w_{p}(t)$ of $\mbox{Ker}\,\mbox{ad}^{r}T_{t}^{lk}\cdot S_{t}^{lk}$, for $t$ from the neighborhood of $0$ in $[0,1)$ and
\begin{gather*}
\mbox{ad}^{r}T^{lk}(t)\cdot S^{lk}(t)\cdot w_{j}(t)=0,\,\,\,j=1,\dots,p.
\end{gather*}
Thus, for $\underline{x}=\varepsilon(t)$ near $\underline{a}$, there exist $w_{1}(\underline{x}),\dots,w_{p}(\underline{x})$ which span the linear space of dimension $\geq s$ and
\begin{gather*}
\mbox{ad}^{r}T^{lk}(\underline{x})\cdot S^{lk}(\underline{x})\cdot w_{j}(\underline{x})=0,\,\,\,j=1,\dots,p.
\end{gather*}
By Lemma 3.3, $\mbox{dim}(\mathcal{R}_{\underline{x}}+\widehat{m}_{\underline{\varphi}(\underline{x})}^{k+1})/\widehat{m}_{\underline{\varphi}(\underline{x})}^{k+1}
=\mbox{dim}\,\mbox{Ker}\,\mbox{ad}^{r}T^{lk}(\underline{x})\cdot S^{lk}(\underline{x})$. Thus for each $\underline{x}\in\tilde{L}\setminus\{\underline{a}\}$ we obtain $\mbox{dim}(\mathcal{R}_{\underline{x}}+\widehat{m}_{\underline{\varphi}(\underline{x})}^{k+1})/\widehat{m}_{\underline{\varphi}(\underline{x})}^{k+1}\geq s$, whence $H_{\underline{x}}\leq H_{\underline{a}}$. This ends the proof.
\end{proof}

\end{proof}

\section{Final remarks}

In this article we proved three implications which consist of proof of Theorem 1. The remaining implications can be proven verbatim in the same way as in \cite{[BM-1]} for analytic case. The reason why the remaining parts of proofs can be repeated in the same way as in subanalytic case is that they rely on properties of analytic functions which are also avaliable for quasianalytic functions. 

It is worth to point out that the class of sets definable in quasianalytic class is broader than class of subanalytic sets. By the results from \cite{[Da]}, in any quasianalytic class which is not analytic but contains analytic functions, there exists a function which is nowhere analytic. Thus there are closed definable sets which are not subanlytic (in classical analytic meaning). The proof of existance of nowhere analytic functions can be found in \cite{[Da]}, however this proof is not constructive and uses Baire theorem. On the other hand, in  \cite{[Ja]}, the author provides a constructive examples of quasianalytic functions which are nowhere analytic( or more generally, which are in some class but nowhere in some lower class).

An interesting problem is to specify some class of closed definable sets which are formally semicoherent (or poses any of the equivalent properties considered in this paper). It is not known so far if even closed quasianalytic sets are formally semicoherent.

\end{document}